%% file: main.tex
\documentclass[a4paper]{article}

\input{packages.tex}
\usepackage[a4paper, total={6in, 8in}]{geometry}
\usepackage{csquotes}
\usepackage[backend=bibtex,style=numeric-comp,giveninits=true,doi=false,isbn=false,url=false,eprint=false,maxbibnames=99]{biblatex}

\input{macros.tex}
\graphicspath{{pictures/}}
\doublefalse

\title{{\papertitle}}
\input{header.tex}

\begin{document}
	\maketitle
	
	\begin{abstract}
		\input{parts_abstract.tex}
	\end{abstract}
	
	\noindent\textbf{Keywords: } \keywordOne, \keywordTwo, \keywordThree, \keywordFour, \keywordFive

	\input{parts_intro.tex}
	\input{parts_methods.tex}

	\input{parts_results.tex}
	\input{parts_discussion.tex}

	\input{acknowledgements.tex}

	\begin{appendices}
		\input{parts_si.tex}	
	\end{appendices}
    \newpage
	\printbibliography

\end{document}

%% file: packages.tex
\usepackage[english]{babel}
\usepackage[utf8]{inputenc}
\usepackage[colorinlistoftodos]{todonotes}
\usepackage{amsmath}
\usepackage{amssymb}
\usepackage{amsthm} 
\usepackage{mathtools}
\usepackage{mathrsfs}
\usepackage{siunitx}
\usepackage{algpseudocode}

\usepackage{subfig}

\usepackage{tikz}

\usepackage{tabularx}
\usepackage{booktabs}
\usepackage{multirow}

\usepackage{empheq}

\usepackage{lipsum}

\usepackage[toc,page]{appendix}

\usepackage{rotating}

\usepackage{makecell}

\usepackage{lscape}
\usepackage{xcolor,colortbl}
\definecolor{lavender}{rgb}{0.9, 0.9, 0.98}

\theoremstyle{plain}
\theoremstyle{plain}
\theoremstyle{plain}
\theoremstyle{plain}
\theoremstyle{plain}
\theoremstyle{definition}
\theoremstyle{remark}

\everymath{\displaystyle}

%% file: macros.tex
\newif\ifdouble

\newcommand{\papertitle}{Real-time whole-heart electromechanical simulations using Latent Neural Ordinary Differential Equations}

\newcommand{\keywordOne}{Cardiac electromechanics}
\newcommand{\keywordTwo}{Machine Learning}
\newcommand{\keywordThree}{Global sensitivity analysis}
\newcommand{\keywordFour}{Parameter estimation}
\newcommand{\keywordFive}{Uncertainty quantification}

\newcommand{\numSims}{N_{\text{sims}}}
\newcommand{\numTrainValid}{N_{\text{train, valid}}}

\newcommand{\numTest}{N_{\text{test}}}

\newcommand{\expected}{\mathbb{E}}
\newcommand{\variance}{\mathbb{V}\mathrm{ar}}

\newcommand{\modFOM}{\mathcal{M}_{\text{3D-0D}}}
\newcommand{\modROM}{\mathcal{M}_{\text{ANN}}}

\newcommand{\ANNState}{\mathbf{z}}
\newcommand{\ANNStateTilde}{\widetilde{\mathbf{z}}}
\newcommand{\ANNLatent}{\mathbf{z}_\mathrm{latent}}

\newcommand{\ANNRhs}{\mathcal{AN \kern-0.3em N}}

\newcommand{\NumANNState} {N_{z}}
\newcommand{\NumANNWeights} {N_{w}}

\newcommand{\qoi}  {\mathbf{q}}

\newcommand{\NoiseCov}{\boldsymbol{\Sigma}}

\newcommand{\ANNparam}{\mathbf{w}}
\newcommand{\ANNparamTrained}{\widehat{\ANNparam}}

\newcommand{\Rtwo}{R\textsuperscript{2}}

\newcommand{\SobolFirst}[2]{S_{#1}^{#2}}
\newcommand{\SobolTotal}[2]{S_{#1}^{#2, T}}

\newcommand{\param}{\boldsymbol{\theta}}

\newcommand{\paramSpace}{\boldsymbol{\Theta}}

\newcommand{\NumSamples} {N_\mathcal{S}}
\newcommand{\NumParams} {N_\mathcal{P}}
\newcommand{\NumMomentum} {N_{\rho}}

\newcommand{\adjoint} {\mathbf{a}}

\newcommand{\momentum}{\boldsymbol{\rho}}

\newcommand{\THB}{T_\mathrm{HB}}

\newcommand{\VLA}{V_{\mathrm{LA}}}
\newcommand{\VLV}{V_{\mathrm{LV}}}
\newcommand{\VRA}{V_{\mathrm{RA}}}
\newcommand{\VRV}{V_{\mathrm{RV}}}

\newcommand{\PLA}{p_{\mathrm{LA}}}
\newcommand{\PLV}{p_{\mathrm{LV}}}
\newcommand{\PRA}{p_{\mathrm{RA}}}
\newcommand{\PRV}{p_{\mathrm{RV}}}

\newcommand{\TLV}{\mathcal{T}_{\text{LV}}}
\newcommand{\Tventricles}{\mathcal{T}_{\text{ventricles}}}
\newcommand{\Tatria}{\mathcal{T}_{\text{atria}}}
\newcommand{\Tall}{\mathcal{T}_{\text{all}}}

\newcommand{\PCabToRORd}{PCa_\mathrm{b}^\mathrm{ToRORd}}
\newcommand{\trpnmaxToRORd}{TRPN^\mathrm{max, ToRORd}}
\newcommand{\GncxbToRORd}{GNCX_\mathrm{b}^\mathrm{ToRORd}}
\newcommand{\TrefToRORdLand}{T_\mathrm{ref}^\mathrm{ToRORd-Land}}
\newcommand{\permfiftyToRORdLand}{perm_\mathrm{50}^\mathrm{ToRORd-Land}}
\newcommand{\npermToRORdLand}{nperm^\mathrm{ToRORd-Land}}
\newcommand{\TRPNnToRORdLand}{TRPN_\mathrm{n}^\mathrm{ToRORd-Land}}
\newcommand{\drToRORdLand}{dr^\mathrm{ToRORd-Land}}
\newcommand{\wfracToRORdLand}{w_\mathrm{frac}^\mathrm{ToRORd-Land}}
\newcommand{\TOTAToRORdLand}{TOT_\mathrm{A}^\mathrm{ToRORd-Land}}
\newcommand{\ktmunblockToRORdLand}{ktm_\mathrm{unblock}^\mathrm{ToRORd-Land}}
\newcommand{\cafiftyToRORdLand}{ca_\mathrm{50}^\mathrm{ToRORd-Land}}
\newcommand{\muToRORdLand}{mu^\mathrm{ToRORd-Land}}
\newcommand{\IupmaxCRN}{I_\mathrm{up}^\mathrm{max, CRN}}
\newcommand{\TrpnmaxCRN}{TRPN^\mathrm{max, CRN}}
\newcommand{\gCaLCRN}{g_\mathrm{CaL}^\mathrm{CRN}}
\newcommand{\TrefCRNLand}{T_\mathrm{ref}^\mathrm{CRN-Land}}
\newcommand{\permfiftyCRNLand}{perm_\mathrm{50}^\mathrm{CRN-Land}}
\newcommand{\npermCRNLand}{nperm^\mathrm{CRN-Land}}
\newcommand{\TRPNnCRNLand}{TRPN_\mathrm{n}^\mathrm{CRN-Land}}
\newcommand{\drCRNLand}{dr^\mathrm{CRN-Land}}
\newcommand{\wfracCRNLand}{w_\mathrm{frac}^\mathrm{CRN-Land}}
\newcommand{\TOTACRNLand}{TOT_\mathrm{A}^\mathrm{CRN-Land}}
\newcommand{\phiCRNLand}{\phi^\mathrm{CRN-Land}}
\newcommand{\cafiftyCRNLand}{ca_\mathrm{50}^\mathrm{CRN-Land}}
\newcommand{\muCRNLand}{mu^\mathrm{CRN-Land}}
\newcommand{\CVventricles}{CV^\mathrm{ventricles}}
\newcommand{\kFEC}{k_\mathrm{FEC}}
\newcommand{\CVatria}{CV^\mathrm{atria}}
\newcommand{\kBB}{k_\mathrm{BB}}
\newcommand{\AVdelay}{AV_\mathrm{delay}}
\newcommand{\Rsys}{R^\mathrm{sys}}
\newcommand{\Rpulm}{R^\mathrm{pulm}}
\newcommand{\Aol}{Aol}
\newcommand{\kArt}{k^\mathrm{Art}}
\newcommand{\aventricles}{a^\mathrm{ventricles}}
\newcommand{\btventricles}{b_\mathrm{t}^\mathrm{ventricles}}
\newcommand{\aatria}{a^\mathrm{atria}}
\newcommand{\bfatria}{b_\mathrm{f}^\mathrm{atria}}
\newcommand{\btatria}{b_\mathrm{t}^\mathrm{atria}}
\newcommand{\kperi}{k_\mathrm{peri}}
\newcommand{\alvrv}{a^\mathrm{lvrv}}
\newcommand{\Treflvrv}{T_\mathrm{ref}^\mathrm{lvrv}}

%% file: header.tex
\author{Matteo Salvador$^{1, 2, 3, *}$,
		Marina Strocchi$^{2, 4}$,
        Francesco Regazzoni$^3$, \\
		Luca Dede'$^3$,
	    Steven A. Niederer$^{2, 4, 5}$,
		Alfio Quarteroni$^{3, 6}$}

\date{\footnotesize
	$^1$ Institute for Computational and Mathematical Engineering, Stanford University, California, USA \\
    $^2$ School of Biomedical Engineering and Imaging Sciences, King's College London, London, UK \\
	$^3$ MOX, Department of Mathematics, Politecnico di Milano, Milan, Italy \\
	$^4$ National Heart and Lung Institute, Imperial College London, London, UK \\
	$^5$ The Alan Turing Institute, London, UK \\
	$^6$ \'Ecole Polytechnique F\'ed\'erale de Lausanne, Lausanne, Switzerland (\textit{Professor Emeritus})\\[2ex]
	$^*$ \textit{Corresponding author} (\texttt{msalvad@stanford.edu}) \\
    }

%% file: parts_abstract.tex
Cardiac digital twins provide a physics and physiology informed framework to deliver predictive and personalized medicine.
However, high-fidelity multi-scale cardiac models remain a barrier to adoption due to their extensive computational costs and the high number of model evaluations needed for patient-specific personalization.
Artificial Intelligence-based methods can make the creation of fast and accurate whole-heart digital twins feasible. 
In this work, we use Latent Neural Ordinary Differential Equations (LNODEs) to learn the temporal pressure-volume dynamics of a heart failure patient.
Our surrogate model based on LNODEs is trained from 400 3D-0D whole-heart closed-loop electromechanical simulations while accounting for 43 model parameters, describing single cell through to whole organ and cardiovascular hemodynamics.
The trained LNODEs provides a compact and efficient representation of the 3D-0D model in a latent space by means of a feedforward fully-connected Artificial Neural Network that retains 3 hidden layers with 13 neurons per layer and allows for 300x real-time numerical simulations of the cardiac function on a single processor of a standard laptop.
This surrogate model is employed to perform global sensitivity analysis and robust parameter estimation with uncertainty quantification in 3 hours of computations, still on a single processor.
We match pressure and volume time traces unseen by the LNODEs during the training phase and we calibrate 4 to 11 model parameters while also providing their posterior distribution.
This paper introduces the most advanced surrogate model of cardiac function available in the literature and opens new important venues for parameter calibration in cardiac digital twins.

%% file: parts_intro.tex
\section{Introduction}
\label{sec:introduction}

Cardiac digital twins integrate physiological and pathological patient-specific data to monitor, analyze and forecast patient disease progression and outcomes.
High-fidelity multi-scale and anatomically accurate models are available but require extensive high-performance computing resources to run, which limit their clinical translation \cite{Niederer2019Nature}.
Over the past years, these mathematical models evolved from an electromechanical description of the human ventricular activity in idealized shapes \cite{landajuela2018numerical, Salvador2020} and realistic geometries \cite{Augustin2021, Regazzoni2022, Piersanti2022, SainteMarie2006}, while also addressing diseased conditions \cite{Salvador2021, Salvador2021MEF, Trayanova2011}, to whole-heart function \cite{Augustin2016, Fedele2023, Gerach2021, Pfaller2019, Peirlinck2021, Strocchi2020Cohort}.
Nevertheless, running many electromechanical simulations still entail high computational costs, hindering the development and application of cardiac digital twins.
The use of Machine Learning tools, such as Gaussian Processes Emulators \cite{Longobardi2020} and Artificial Neural Networks (ANNs) \cite{regazzoni2019modellearning, Salvador2021KPOD}, allows to create efficient surrogate models that can be employed in many-query applications \cite{QMN16}, such as sensitivity analysis and parameter inference \cite{Regazzoni2022EMROM, Salvador2023, Strocchi2023}.
In the framework of digital twinning and personalized medicine, bridging the chasm between the need for a supercomputer \cite{Marchesseau2013, Marx2020, Sermesant2012, Strocchi2020Simulating} and performing accurate real-time numerical simulations on a standard computer \cite{Cicci2022, Jung2022, Regazzoni2021cardioemulator, Schiavazzi2017} would have a tremendous impact on the clinical adoption of computational cardiology.

In this work, we develop a Scientific Machine Learning method to build, to the best of our knowledge, the most comprehensive surrogate model involving both cardiac and cardiovascular function that is currently available in the literature.
Specifically, we train a system of Latent Neural Ordinary Differential Equations (LNODEs) \cite{Chen2019, Rubanova2019, regazzoni2019modellearning} that learns the pressure-volume transients of a heart failure patient while varying 43 model parameters that describe cardiac electrophysiology, active and passive mechanics, and cardiovascular fluid dynamics, by employing 400 3D-0D closed-loop electromechanical training simulations.
We design a suitable loss function that is minimized during the tuning process of the ANN parameters, which entails small relative errors of LNODEs, i.e. from $2\%$ to $5\%$, when the number of training samples is small compared to the dimensionality of the parameter space and the explored model variability.
These LNODEs allow for 300x real-time four-chamber heart numerical simulations and can be easily trained on a single Central Processing Unit (CPU).

We use the trained LNODEs to perform global sensitivity analysis (GSA) and robust parameter estimation with uncertainty quantification (UQ) \cite{Regazzoni2022EMROM, Salvador2023}.
For the former, we observe how model parameters impact the variability of scalar quantities of interest (QoIs) retrieved from the pressure-volume time traces, by considering both first-order and high-order interactions via Sobol indices \cite{sobol1990sensitivity}.
For the latter, we combine two Bayesian statistics methods, i.e. Maximum a Posteriori (MAP) estimation and Hamiltonian Monte Carlo (HMC) \cite{Chen2019, Betancourt2017, Homan2014}, where we exploit efficient matrix-free adjoint-based methods, automatic differentiation and vectorization \cite{Chen2019}.
In particular, we design several test cases where we calibrate tens of model parameters by matching the pressure and volume time traces, that are time-dependent QoIs, coming from 5 unseen 3D-0D numerical simulations for the trained ANN.
GSA and parameter estimation with UQ can be carried out in 3 hours of computations by using a single core standard laptop.

%% file: parts_methods.tex
\section{Methods}
\label{sec:methods}

\begin{figure}[t!]
    \centering
    \includegraphics[width=1.0\textwidth]{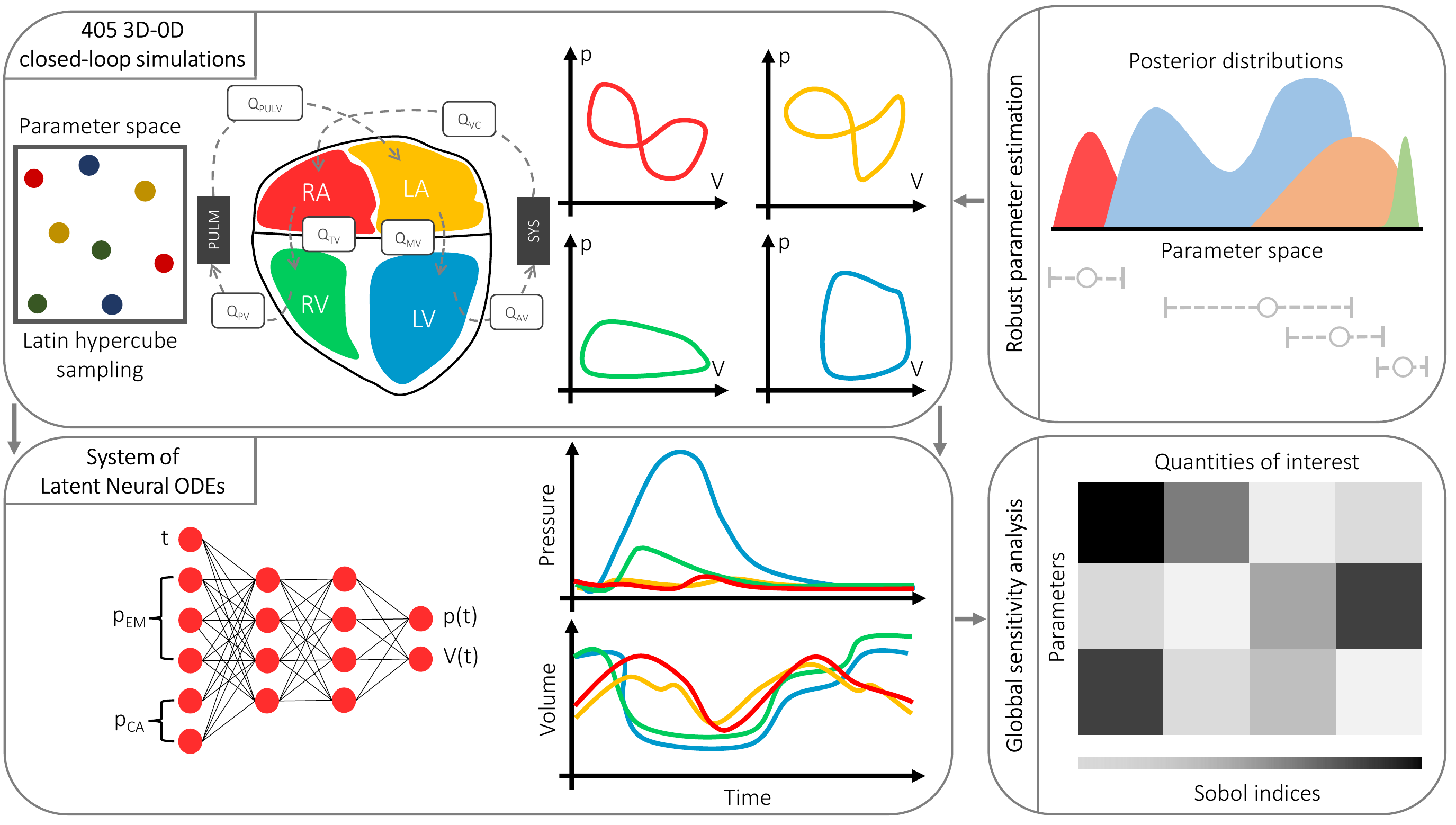}
    \caption{Sketch of the computational pipeline. We perform several 3D-0D closed-loop four-chamber heart electromechanical simulations. We build an accurate and efficient ANN-based surrogate model of the whole cardiovascular function by means of LNODEs. We carry out GSA to understand how each model parameter influences different QoIs extracted from the simulated pressure-volume loops. We robustly estimate many model parameters from time-dependent QoIs. Fully personalized 3D-0D numerical simulations can be performed after parameters calibration with UQ.}
    \label{fig:computationalpipeline}
\end{figure}

We display the whole computational pipeline in Figure~\ref{fig:computationalpipeline}. 
\begin{itemize}
\item Top-left: we use a database of $\numSims = 405$ electromechanical simulations generated by a personalized anatomy four-chamber heart model from a heart failure patient (see Appendix~\ref{app:geometry}), where we vary $\NumParams = 43$ parameters that describe cell, tissue, whole-heart and cardiovascular system material properties and boundary conditions. For all the numerical simulations, we run 5 heartbeats in sinus rhythm and we perform our analysis on the pressure and volume transients of the last cardiac cycle. We refer to Appendix~\ref{app:methods:3D0D} for all the details about the four-chamber physics-based mathematical model and the numerical settings of these simulations. All the information regarding model parameters can be found in Appendix~\ref{app:methods:parameters}.
\item Bottom-left: we employ $\numTrainValid = 400$ simulations to tune the LNODEs hyperparameters. This surrogate model learns the atrial and ventricular pressure-volume temporal dynamics of the last cardiac cycle only, while receiving time and model parameters as inputs. We perform $K$-fold cross validation with $K = 10$ for the training-validation splitting. We detail the whole optimization process to get the final values of the LNODEs hyperparameters in Appendix~\ref{app:methods:LNODE}. We evaluate the accuracy of the trained LNODEs on a testing dataset consisting of the remaining $\numTest = 5$ numerical simulations.
\item Bottom-right: we employ the trained LNODEs to perform GSA.
\item Top-right: we estimate model parameters with UQ on $\numTest = 5$ numerical simulations by means of the trained LNODEs.
\end{itemize}

\subsection{Learning atrial and ventricular pressure-volume loops}
\label{sec:methods:LNODE}

Following the model learning approach introduced in \cite{regazzoni2019modellearning}, we build a system of LNODEs, i.e. a set of ordinary differential equations whose right hand side is represented by a feedforward fully-connected ANN, that learns the pressure-volume temporal dynamics of the 3D-0D closed-loop electromechanical model $\modFOM$ in a latent space.
In this framework, the four-chamber heart surrogate model $\modROM$ reads: 
\begin{equation} \label{eqn:ROM}
    \left\{
    \begin{aligned}
        \frac{d \ANNState(t)}{d t} &= \ANNRhs\left(
            \ANNState(t),
            \cos\left(\tfrac{2 \pi (t - \AVdelay)}{\THB}\right),
            \sin\left(\tfrac{2 \pi (t - \AVdelay)}{\THB}\right),
            \param;
            \ANNparam
        \right)
        && \text{for } t \in (0, \THB],\\
        \ANNState(0) &= \ANNState_{0},  && \\
    \end{aligned}
    \right.
\end{equation}
where $\ANNState_{0}$ is the vector of initial conditions.
The ANN, with weights and biases encoded in $\ANNparam \in \mathbb{R}^{\NumANNWeights}$, is defined by $\ANNRhs \colon \mathbb{R}^{\NumANNState + 2 + \NumParams} \to \mathbb{R}^{\NumANNState}$.
Vector $\param \in \paramSpace \subset \mathbb{R}^{\NumParams}$ defines the model $\modFOM$ parameters. Some examples of $\param$ could be conductances of different ionic channels, myocardial conductivity, atrial and ventricular active tension or passive stiffness, and resistances of the systemic and pulmonary circulation.
The reduced state vector $\ANNState(t) \in \mathbb{R}^{\NumANNState}$ contains the time-dependent pressure and volume variables of the left atrium (LA), right atrium (RA), left ventricle (LV) and right ventricle (RV), as well as additional latent variables without a direct physical interpretation, that is $\ANNState(t) = [\PLA(t), \PLV(t), \PRA(t), \PRV(t), \VLA(t), \VLV(t), \VRA(t), \VRV(t), \ANNLatent(t)]^T$.
The ANN receives $\NumANNState$ state variables, $\NumParams$ scalar parameters, and two periodic inputs.
Indeed, even though LNODEs are just trained on the last cardiac cycle, the $\cos({2 \pi (t - \AVdelay)}/{\THB})$ and $\sin({2 \pi (t - \AVdelay)}/{\THB})$ terms account for the heartbeat period $\THB$ and the atrioventricular delay $\AVdelay$ of whole-heart electromechanical simulations (see Appendix~\ref{app:methods:3D0D} for further details).
We stress that, differently from \cite{Regazzoni2022EMROM}, the initial reduced state vector $\ANNState_0$ contains different sets of initial conditions for pressures, volumes and latent variables \cite{Regazzoni2021assimilation}. 

The loss function that we minimize during the ANN optimization process reads:
\begin{equation} \label{eqn:loss}
\begin{split}
\mathcal{L}(\ANNState(t), \ANNStateTilde(t); \ANNparamTrained) = \underset{\ANNparamTrained}{\arg\min} &\biggl[ \dfrac{|| \ANNState(t) - \ANNStateTilde(t) ||_{\text{L}^2(0, \THB)}^2}{\ANNState_{\text{norm}}^2} \\
&+ \alpha \dfrac{\displaystyle \Big\lvert \Big\lvert \frac{d \ANNState(t)}{d t} - \frac{d \ANNStateTilde(t)}{d t} \Big\lvert \Big\lvert_{\text{L}^2(0, \THB)}^2}{\ANNState_{\text{norm, diff}}^2} \\
&+ \beta \dfrac{\left(\underset{t \in [0, \THB]}{\max} \ANNState(t) - \underset{t \in [0, \THB]}{\max} \ANNStateTilde(t)\right)^2}{\ANNState_{\text{norm, max}}^2} \\
&+ \gamma \dfrac{\left(\underset{t \in [0, \THB]}{\min} \ANNState(t) - \underset{t \in [0, \THB]}{\min} \ANNStateTilde(t)\right)^2}{\ANNState_{\text{norm, min}}^2} \\
&+ \eta \left( || \ANNLatent(0) ||^2 + || \ANNLatent(\THB) ||^2 \right) \\
&+ \iota || \ANNparamTrained ||_{\text{L}^2}^2 \biggr],
\end{split}
\end{equation}
with $\alpha = \beta = \gamma = \eta = 0.1$.
The loss function aims at finding an optimal set of weights $\ANNparamTrained$ for the ANN.
It comprises the normalized mean square error between ANN predictions $\ANNState(t)$ and observations $\ANNStateTilde(t)$, as well as a weak penalization of the reduced state vector time derivatives, maximum and minimum values for $t \in [T - \THB, T]$.
Indeed, given the small ratio between the dimensionality of the training dataset and the number of parameters $\param$ of model $\modFOM$, we notice that these three additional terms reduce the generalization errors of the ANN.
The penultimate weakly enforced condition on $\ANNLatent(t)$ favors a periodic solution for all the hidden latent variables.
The last term of the loss function prescribes the $L^2$ regularization of the ANN weights and $\iota$ is one of the automatically tuned LNODEs hyperparameters (see Appendix~\ref{app:methods:LNODE}).

\subsection{Global sensitivity analysis}
\label{sec:methods:GSA}
We employ the Saltelli's method to perform a variance-based sensitivity analysis \cite{saltelli2002making}.
We compute both first-order Sobol indices and total-effect Sobol indices for each combination of quantity of interest and model parameter \cite{sobol1990sensitivity}.
These two indices define how much varying a single parameter affects a specific QoI and how higher-order interactions among model parameters influences the model outputs, respectively.
Further details are provided in Appendix~\ref{app:methods:GSA}.

\subsection{Robust parameter estimation}
\label{sec:methods:PE}
We perform parameter calibration with inverse UQ following a two-stage approach.
First, given a set of time-dependent QoIs related to four-chamber heart pressure and volume traces, we solve a bounded and constrained optimization problem by employing model $\modROM$ to obtain the pointwise MAP estimation for a predefined set of model parameters $\param \in \paramSpace \subset \mathbb{R}^{\NumParams}$.
Second, we initialize HMC based on the MAP estimation and we build an approximation for the posterior distribution of $\param$ \cite{Betancourt2017}, while accounting for the measurement and surrogate modeling errors via Gaussian Processes \cite{Salvador2023}.
We provide all the mathematical and numerical details about these techniques in Appendix~\ref{app:methods:PE}.

\subsection{Software libraries}
\label{sec:methods:software}
All 3D-0D closed-loop electromechanical simulations run with the Cardiac Arrhythmia Research Package (CARP) \cite{Augustin2016, Vigmond2003}.
We train model $\modROM$ by using an in-house high-performance Python library based on Tensorflow \cite{Tensorflow2015}.
We perform GSA by means of the open source Python library \texttt{SALib}\footnote{\url{https://salib.readthedocs.io/}} \cite{Herman2017}.
Parameter estimation with UQ is carried out by combining the open source Python libraries \texttt{JAX}\footnote{\url{https://github.com/google/jax}} \cite{jax2018github} and \texttt{NumPyro}\footnote{\url{https://github.com/pyro-ppl/numpyro}} \cite{phan2019}.
This paper is accompanied by \url{https://github.com/MatteoSalvador/cardioEM-4CH}, a public repository containing the trained LNODEs, along with the codes to perform GSA and robust parameter identification.

%% file: parts_results.tex
\section{Results}
\label{sec:results}

We provide the numerical results for the training and testing phases of LNODEs, along with their application to GSA and robust parameter estimation.

\subsection{Learning atrial and ventricular pressure-volume loops}
\label{sec:results:LNODE}

\begin{figure}
    \centering
    \includegraphics[width=1.0\textwidth]{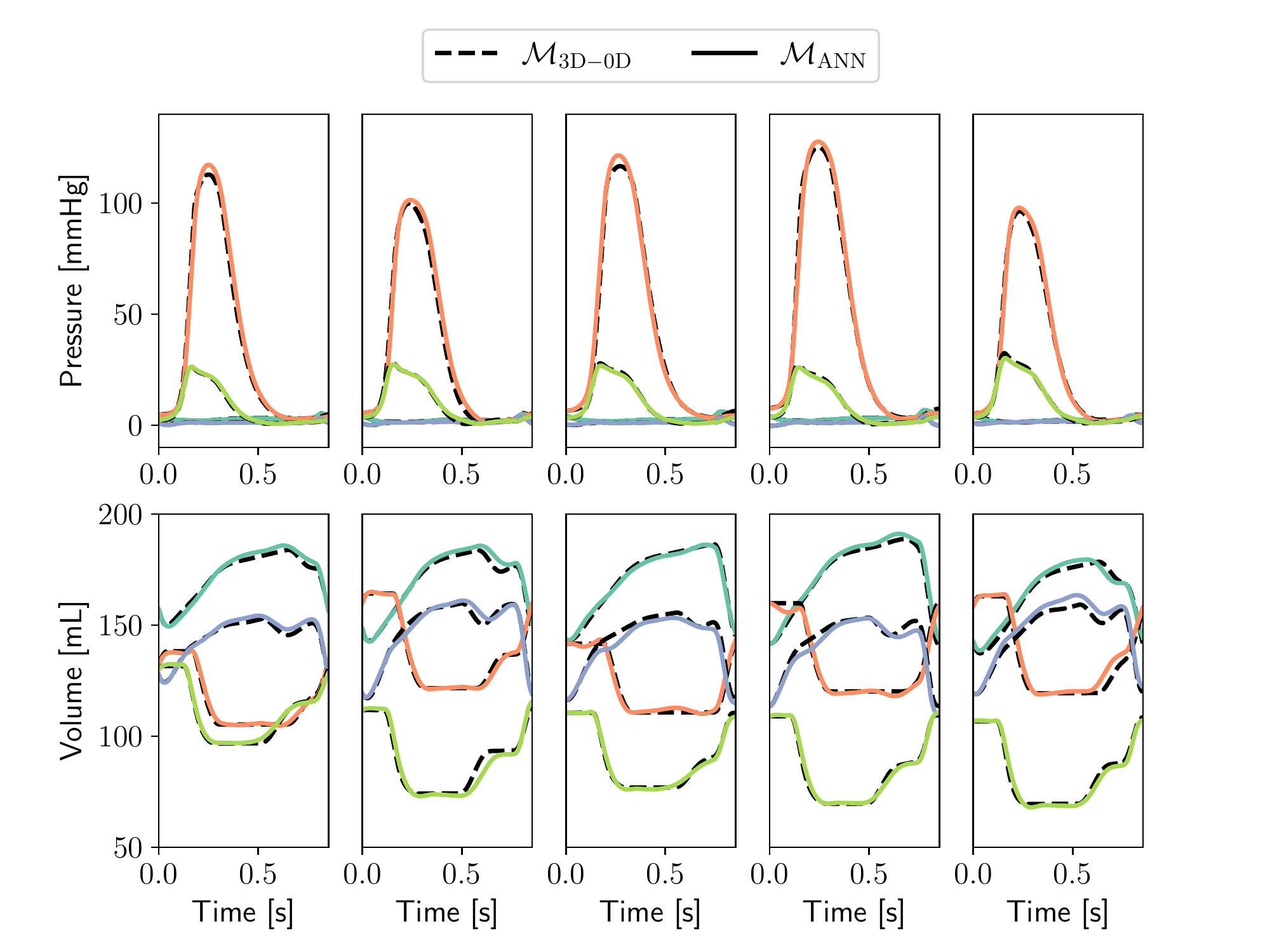}
    \caption{Pressure and volume time transients obtained with $\modFOM$ (dashed lines), compared to those obtained with model $\modROM$ (solid lines), on the testing samples ($\numTest = 5$). Light blue: LA, orange: LV, blue: RA, green: RV.}
    \label{fig:pressurevolume}
\end{figure}

\begin{table}
\begin{center}
       \begin{tabular}{ rlrrrr }
       \toprule
       & & \multicolumn{4}{c}{\textbf{Pressure}} \\
       & & $\PLA(t)$ & $\PLV(t)$ & $\PRA(t)$ & $\PRV(t)$
       \\
       \multirow{ 2}{*}{$\modFOM$ vs $\modROM$} & \multicolumn{1}{l|}{NRMSE} &
       0.027522 & 0.021890 & 0.021776 & 0.020445
       \\
       & \multicolumn{1}{l|}{\Rtwo} & 99.2319 & 99.8189 & 98.8457 & 99.8139
       \\
       \midrule
       & & \multicolumn{4}{c}{\textbf{Volume}} \\
       & & $\VLA(t)$ & $\VLV(t)$ & $\VRA(t)$ & $\VRV(t)$
       \\
       \multirow{ 2}{*}{$\modFOM$ vs $\modROM$} & \multicolumn{1}{l|}{NRMSE} & 0.035943 & 0.030143 & 0.054243 & 0.026157
       \\
       & \multicolumn{1}{l|}{\Rtwo} & 99.3619 & 99.4978 & 97.9704 & 99.5761
       \\
       \bottomrule
       \end{tabular}
       \caption{Testing errors and \Rtwo coefficients on the time-dependent outputs of the trained LNODEs system.}
       \label{tab:errorstime}
\end{center}
\end{table}

Automatic hyperparameters tuning with $K$-fold cross validation leads to an optimal ANN architecture comprising 3 hidden layers and 13 neurons per hidden layer.
The optimal number of states is set to $\NumANNState = 8$, i.e. no latent variables are selected.
This is motivated by the trade-off between the size of the training set $\numTrainValid$ with respect to the number of parameters $\NumParams$, i.e. a thrifty system of LNODEs with no additional hidden variables $\ANNLatent(t)$ is selected to avoid overfitting.
More details regarding LNODEs training and hyperparameters tuning are given in Appendix~\ref{app:methods:LNODE}.

In Table~\ref{tab:errorstime}, we report the Normalized Root Mean Square Error (NRMSE) and $\Rtwo$ coefficients associated with the LA, LV, RA and RV pressure-volume time traces provided by LNODEs.
These values are obtained by considering a test set comprised of $\numTest = 5$ electromechanical simulations.
The accuracy obtained by our surrogate model in reproducing the cardiac outputs is high, manifesting testing errors that approximately range from $2\%$ to $5\%$ for all time-dependent QoIs.
The good match between models $\modFOM$ and $\modROM$ is also confirmed by Figure~\ref{fig:pressurevolume}, where atrial and ventricular pressure-volume traces present a good overlap on the whole testing set.

\subsection{Global sensitivity analysis}
\label{sec:GSA}

\begin{figure}
    \centering
    \includegraphics[width=1.1\textwidth]{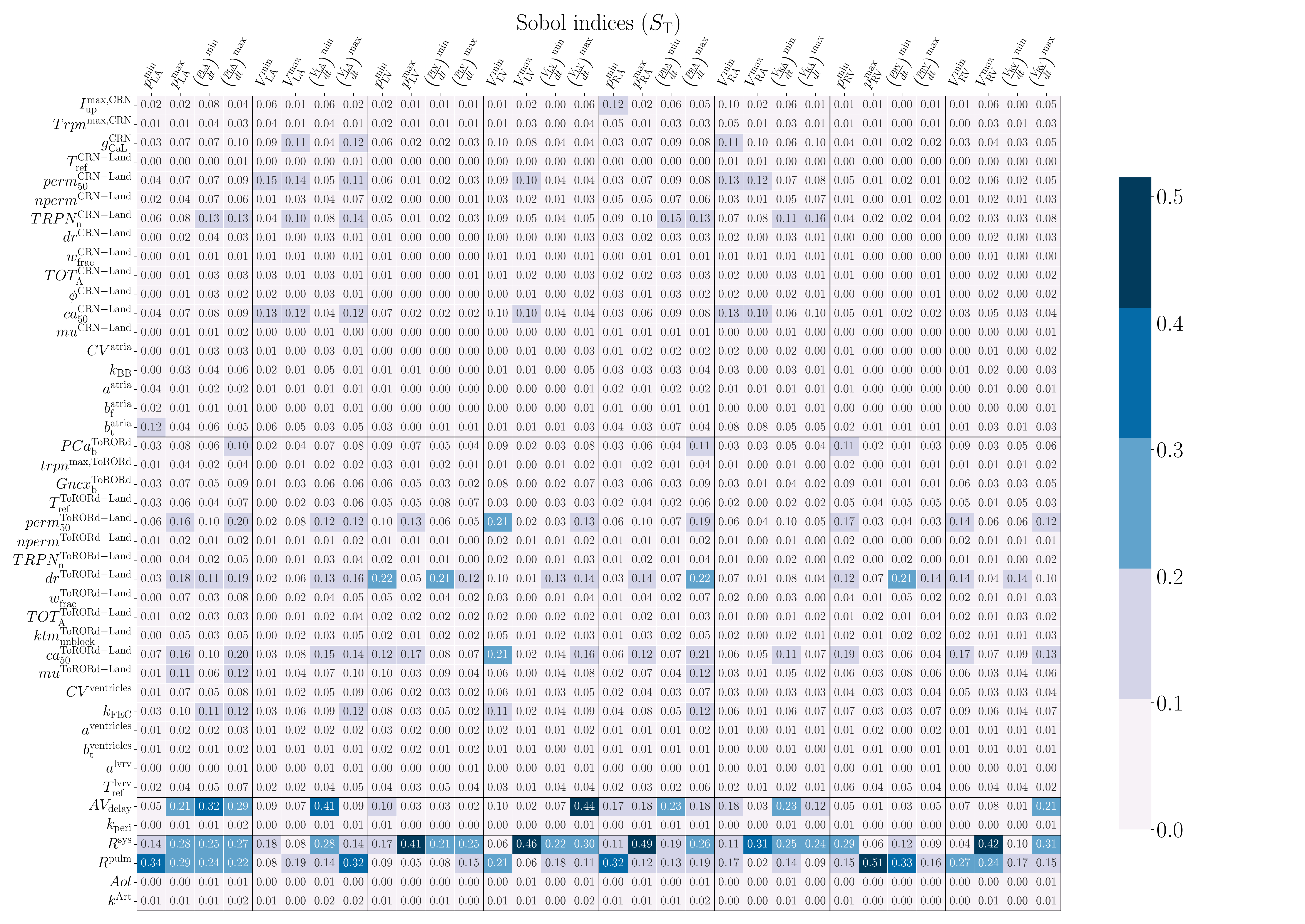}
    \caption{Total-effect Sobol indices computed by exploiting model $\modROM$.
             For a detailed definition of all model parameters and QoIs, we refer to Appendix~\ref{app:methods:parameters}.}
    \label{fig:ST}
\end{figure}

Figure~\ref{fig:ST} shows the total-effect Sobol indices.
We consider a parameter to be relevant if the associated Sobol indices are greater than $10^{-1}$ for at least one QoI.
We notice that some model parameters are compartmentalized, i.e. cell-to-organ level values coming from a certain compartment of the cardiocirculatory system mostly explain the variability of QoIs that are specific to that region.
Indeed, some parameters of the CRN-Land model, such as $\permfiftyCRNLand$, $\TRPNnCRNLand$ and $\cafiftyCRNLand$, or of the Guccione model, such as $\btatria$, have an important role in determining atrial behavior.
Similar considerations occur for the ventricular part of the heart, where the most important parameters are related to the ToRORd-Land model.
Nevertheless, it is important to notice the interplay between some ventricular parameters of the ToRORd-Land model at the cellular scale, such as $\drToRORdLand$, $\permfiftyToRORdLand$ and $\cafiftyToRORdLand$ and the atrial function.
Finally, we highlight that some model parameters, such as atrioventricular delay $\AVdelay$, systemic resistance $\Rsys$ and pulmonary resistance $\Rpulm$ strongly affect all QoIs, whereas others, such as the pericardial coefficient $\kperi$, as well as aorta parameters ($\Aol$, $\kArt$), have a minor role in determining all QoIs.

\subsection{Robust parameter estimation}
\label{sec:results:PE}

\begin{table}[t!]
    \begin{center}
        \begin{tabular}{ l|cc }
            \toprule
            Test case & Time-dependent QoIs & Estimated model parameters \\
            \midrule
            $\TLV$          & $\VLV$ & \makecell{$\drToRORdLand$, $\cafiftyToRORdLand$, $\Rsys$, $\Rpulm$} \\ 
            \midrule
            $\Tventricles$  & $\VLV$, $\VRV$ & \makecell{$\drToRORdLand$, $\cafiftyToRORdLand$, \\ $\permfiftyToRORdLand$, $\Rsys$, $\Rpulm$} \\
            \midrule
            $\Tatria$       & $\VLA$, $\VRA$ & \makecell{$\cafiftyCRNLand$, $\permfiftyCRNLand$, $\TRPNnCRNLand$, \\ $\gCaLCRN$, $\Rsys$, $\Rpulm$} \\
            \midrule
            $\Tall$         & \makecell{$\PLA$, $\PRA$, $\PLV$, $\PRV$, \\ $\VLA$, $\VRA$, $\VLV$, $\VRV$} & \makecell{$\drToRORdLand$, $\permfiftyToRORdLand$, $\cafiftyToRORdLand$, \\ $\cafiftyCRNLand$, $\permfiftyCRNLand$, $\TRPNnCRNLand$, \\ $\kFEC$, $\gCaLCRN$, $\btatria$, $\Rsys$, $\Rpulm$} \\
            \bottomrule
        \end{tabular}
        \caption{Summary of the 4 in silico test cases for parameter calibration.}
        \label{tab:testcases}
    \end{center}
\end{table}

In the context of parameter calibration, a preliminary GSA allows to determine the identifiability of model parameters according to the provided QoIs.
Based on the results obtained in Section~\ref{sec:GSA}, we design 4 in silico test cases to show the robustness and flexibility of our parameter calibration process, which is driven by a combined use of MAP estimation and HMC starting from time-dependent QoIs.
In Table~\ref{tab:testcases}, we report the observed pressure-volume time traces and estimated model parameters for each test case.
In $\TLV$ and $\Tventricles$, we estimate model parameters related to the ventricular and cardiovascular function starting from time-dependent QoIs localized in the ventricles.
In $\Tatria$, we calibrate model parameters over the whole cardiac function and cardiocirculatory network by only considering atrial observations.
Finally, we challenge our surrogate model by taking all cardiac pressures and volumes over time and by estimating 11 model parameters.

\begin{figure}
    \centering
    \includegraphics[width=1.0\textwidth]{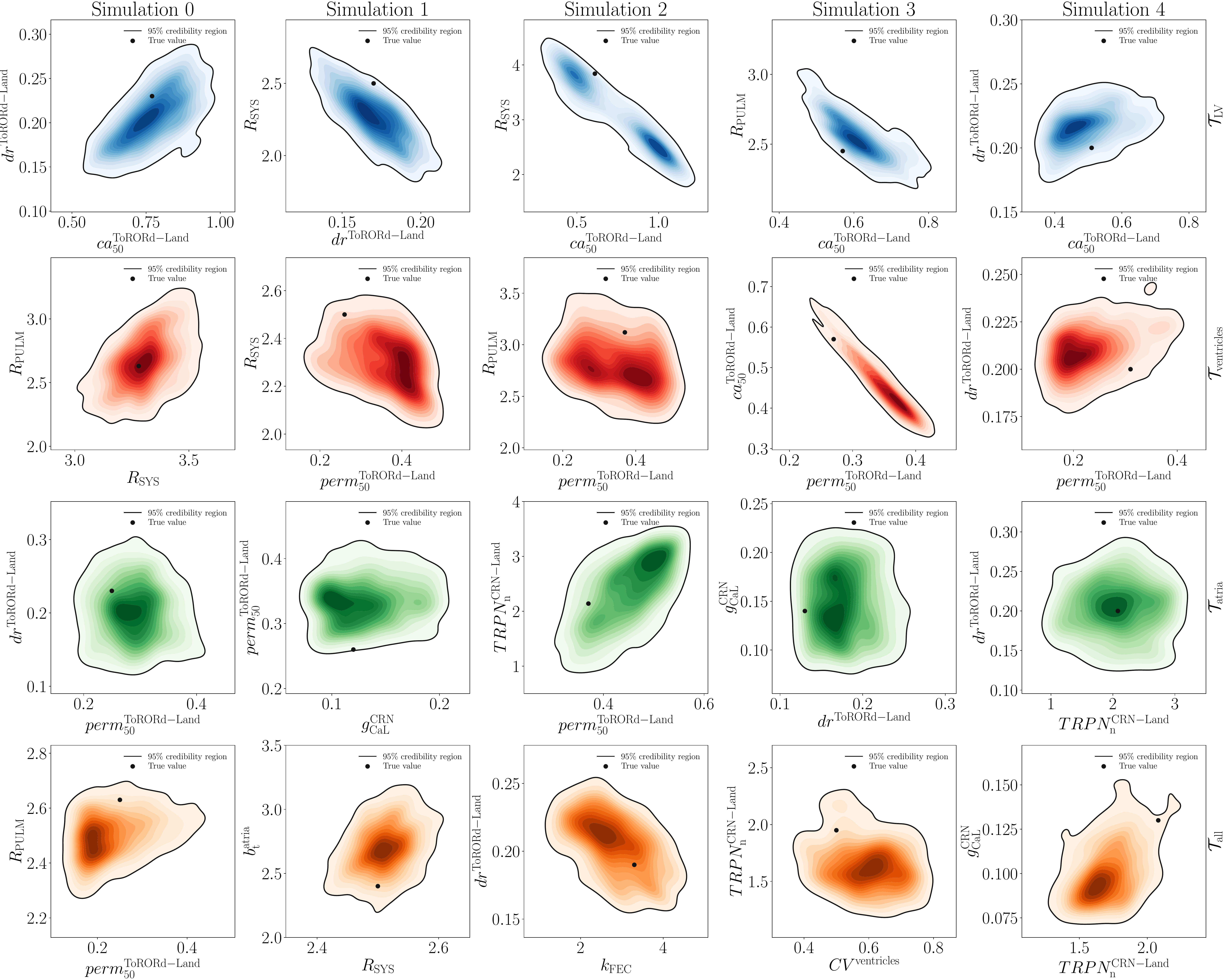}
    \caption{Two-dimensional views of the posterior distribution estimated by means of HMC for each test case (rows) and $\numTest = 5$ electromechanical simulations (columns). Different colors are associated to $\TLV$, $\Tventricles$, $\Tatria$, $\Tall$.}
    \label{fig:posterior}
\end{figure}

We perform parameter estimation with UQ on $\numTest = 5$ electromechanical simulations that are unseen by the trained LNODEs.
Figure~\ref{fig:posterior} shows some two-dimensional views of the posterior distribution for each test case and for all $\numTest$ numerical simulations.
We notice that the true parameter values are contained inside the $95\%$ credibility regions.
Moreover, by using Bayesian statistics we are able to capture relationships among model parameters.
In particular, in Figure~\ref{fig:posterior} we consider different pairs of model parameters for each test case and numerical simulation to maximize the number of interactions.
For instance, $\Rpulm$ and $\drToRORdLand$ are positively correlated with $\Rsys$ and $\cafiftyToRORdLand$, respectively, while $\drToRORdLand$ and $\cafiftyToRORdLand$ are negatively correlated with $\kFEC$ and $\permfiftyToRORdLand$, respectively. 
We notice that, in some cases, cell-based atrial and ventricular parameters may be correlated, as it happens for $\TRPNnCRNLand$ and $\permfiftyToRORdLand$, while in most situations, such as with $\drToRORdLand$ and $\TRPNnCRNLand$, there is no interaction.
We also remark that this kind of relationships may be unraveled among different physical problems.
For instance, this occurs between cardiovascular hemodynamics ($\Rsys$) and the ventricular cell tension model ($\drToRORdLand$).
For the sake of completeness, in Table~\ref{tab:HMC} we report the identified parameter values of $\drToRORdLand$, $\Rsys$ and $\Rpulm$ for all test cases, with respect to the first testing simulation.
We show that the true values of the parameters are always contained inside the interval defined by mean plus/minus two standard deviations.
We refer to Appendix~\ref{app:methods:PE} for the tables containing similar results and comparisons for all test cases ($\TLV$, $\Tventricles$, $\Tatria$ and $\Tall$) with all the relevant model parameters over the $\numTest$ electromechanical simulations.

\begin{table}[t!]
    \begin{center}
        \begin{tabular}{ l|rrrrr }
            \toprule
            Parameter & Ground truth & $\TLV$ & $\Tventricles$ & $\Tatria$ & $\Tall$ \\
            \midrule
            $\drToRORdLand$ & 0.23 & 0.21 $\pm$ 0.06 & 0.23 $\pm$ 0.05 & 0.20 $\pm$ 0.07 & 0.27 $\pm$ 0.05 \\ 
            $\Rsys$         & 3.28 & 3.28 $\pm$ 0.63 & 3.30 $\pm$ 0.20 & 3.33 $\pm$ 0.38 & 3.18 $\pm$ 0.10 \\
            $\Rpulm$        & 2.63 & 2.75 $\pm$ 0.63 & 2.67 $\pm$ 0.41 & 2.98 $\pm$ 0.67 & 2.50 $\pm$ 0.16 \\
            \bottomrule
        \end{tabular}
        \caption{Mean plus/minus two standard deviations associated to the estimated values of $\drToRORdLand$, $\Rsys$ and $\Rpulm$ during HMC for the first sample of the testing set.}
        \label{tab:HMC}
    \end{center}
\end{table}

%% file: parts_discussion.tex
\section{Discussion}
\label{sec:discussion}

\begin{table}[t!]
    \begin{center}
        \hspace*{-0.6cm}
        \begin{tabular}{ l c cc cc c}
            \toprule
            Task & Computational resources & Execution time \\
            \midrule
          	$\modFOM$ \\
            Single simulation (5 heartbeats)               & 512 cores & 6 hours and 20 minutes \\
            GSA (704'000 simulations)                      & 512 cores & 508 years \\
            Parameter estimation with UQ (750 heartbeats)  & 512 cores & 0.5 years \\            
                                                           &           & Total: 508.5 years \\
            \midrule
            $\modROM$ \\
            Training dataset generation (405 simulations)  & 512 cores & 106 days and 21 hours \\
            Reduced-order model training                   & 1 core    & 10 hours \\
            GSA (704'000 heartbeats)                       & 1 core    & 2 hours \\
            Parameter estimation with UQ (750 heartbeats)  & 1 core    & 1 hour \\
                                                           &           & Total: 108 days \\
            \bottomrule
        \end{tabular}
        \caption{Summary of the approximated computational times to perform GSA and parameter estimation with UQ. 3D-0D closed-loop model $\modFOM$ (top) and LNODEs $\modROM$ (bottom).}
        \label{tab:computationaltimes}
    \end{center}
\end{table}

In this work, we propose a surrogate model based on LNODEs to learn the pressure-volume temporal dynamics of 3D-0D closed-loop four-chamber heart electromechanical simulations \cite{Regazzoni2022EMROM}.
Starting from 400 numerical simulations, we create a surrogate model of a heart failure patient by leveraging LNODEs.
These are defined by a lightweight feedforward fully-connected ANN containing 3 hidden layers and 13 neurons per layer.
LNODEs retain the variability of 43 model parameters that describe electrophysiology, active and passive mechanics, and hemodynamics, both at the cell level and organ scale, and covering a wide range of pressure and volume values (see Figures~\ref{fig:pressuresvolumes} and ~\ref{fig:pvloops} in Appendix~\ref{app:methods:3D0D}).
The generation of such a comprehensive training dataset poses an incredible technological challenge itself in the scientific community \cite{Strocchi2023}.
On top of that, this paper provides, to the best of our knowledge, the most comprehensive surrogate model embracing cardiac and cardiovascular function that has been currently proposed in the literature.
With respect to other Machine Learning tools, such as Gaussian Processes Emulators \cite{Longobardi2020}, LNODEs present a higher representational power, because they encode time dependent numerical simulations instead of pointwise QoIs, while also requiring a smaller amount of data to reach a prescribed accuracy \cite{Regazzoni2022EMROM}.

LNODEs require a small amount of computational resources and enable several applications of interest in a very fast and accurate manner.
Indeed, as reported in Table~\ref{tab:computationaltimes}, running the training phase of the ANN along with GSA and robust parameter estimation on a single core standard laptop just requires 13 hours of computations.
We remark that this time can be reduced with a multi-core implementation.
On the other hand, employing the 3D-0D model $\modFOM$ for the same computational pipeline would entail very significant costs.
The overall speed-up with the surrogate model $\modROM$ is equal to 1718x.
The extension of the proposed method to incorporate different anatomies and pathological conditions would potentially allow for a universal whole-heart simulator that might be readily deployed in clinical practice for fast and reliable computational analysis.


%% file: acknowledgements.tex
\section*{Acknowledgements}
This project has been funded by the Italian Ministry of University and Research (MIUR) within the PRIN (Research projects of relevant national interest 2017 “Modeling the heart across the scales: from cardiac cells to the whole organ” Grant Registration number 2017AXL54F).
This project has also been supported by the INdAM GNCS Project CUP E55F22000270001.
SAN acknowledges NIH R01-HL152256, ERC PREDICT-HF 453 (864055), BHF (RG/20/4/34803), EPSRC (EP/P01268X/1, EP/X012603/1), EPSRC Grant EP/X03870X/1 and The Alan Turing Institute.
LD acknowledges the support by the FAIR (Future Artificial Intelligence Research) project, funded by the NextGenerationEU program within the PNRR-PE-AI scheme (M4C2, investment 1.3, line on Artificial Intelligence), Italy.

%% file: parts_si.tex
\section{Four-chamber heart geometry}
\label{app:geometry}

\begin{figure}[t!]
\begin{center}
\includegraphics[width=0.7\textwidth]{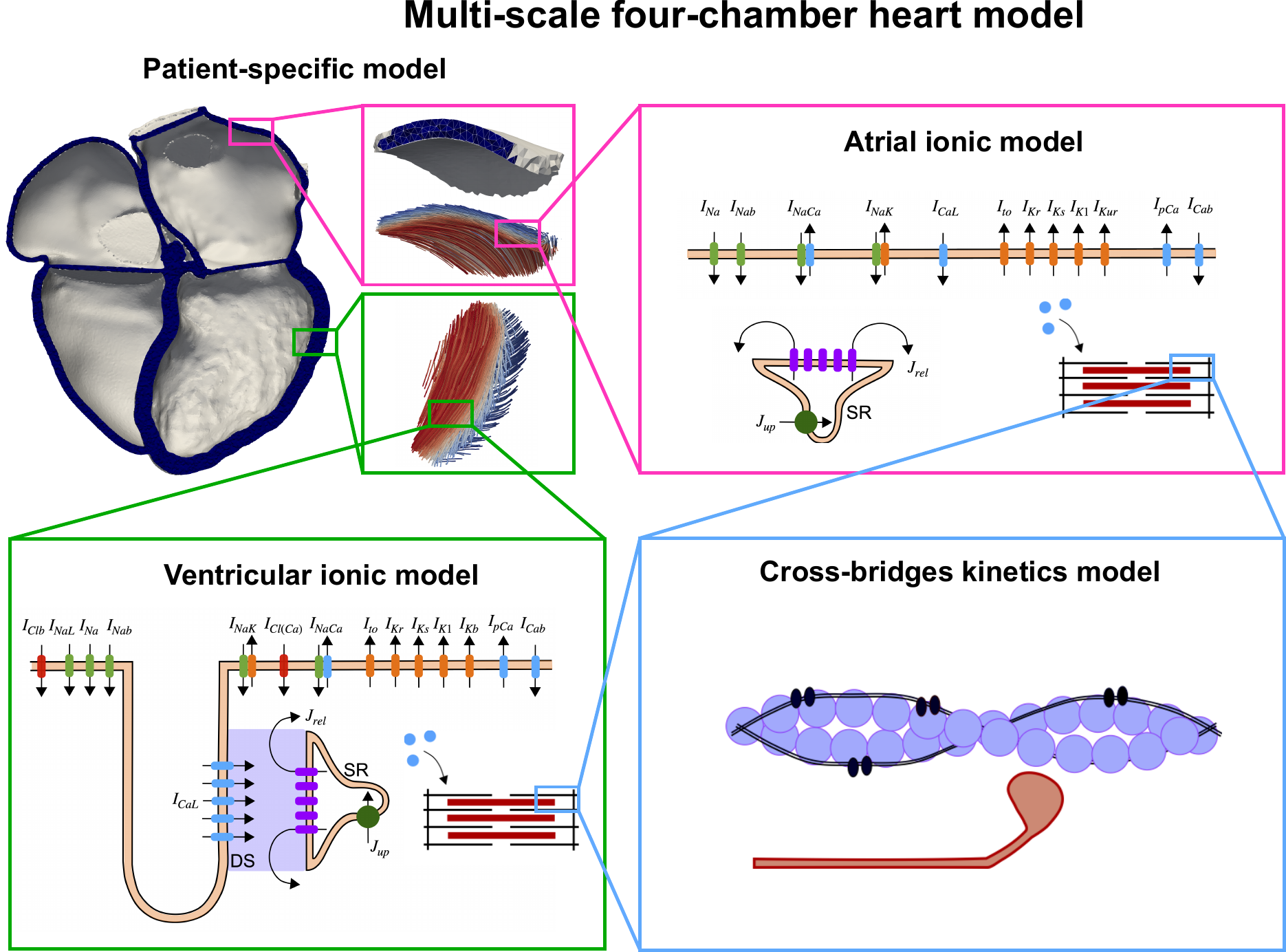}
\end{center}
\caption{Four-chamber heart geometry. Patient-specific whole-heart mesh (left). Refined atria, atrial and ventricular transmural myofiber orientation from endocardium to epicardium (right) \cite{Strocchi2023}.}
\label{fig:geometry}
\end{figure}

The end-diastolic computed tomography (CT) image acquired from a 77 yo female heart failure patient with atrial fibrillation was segmented to generate a four-chamber heart geometry.
All the computational tools regarding segmentation and meshing with 1 mm linear tetrahedral Finite Elements are described in \cite{Strocchi2020Cohort,Strocchi2020Simulating}.
The atria are refined with the resample algorithm from meshtool \cite{Neic2020} to have at least 3 elements across the wall thickness to reduce locking effects.
The ventricles were assigned with a transmural fibre distribution using the Bayer's rule-based algorithm \cite{Bayer2012} (Figure~\ref{fig:geometry}, bottom right), where the fibre and sheet angles at the endocardium and epicardium are +60$^\circ$ and -60$^\circ$ \cite{Pfaller2019}, and -65$^\circ$ and +25$^\circ$ \cite{Bayer2012}, respectively.
Atrial myofibre orientation was assigned by computing universal atrial coordinates on the atria and by mapping an ex-vivo diffusion tensor MRI dataset onto the endocardial and the epicardial surfaces (Figure~\ref{fig:geometry}, top right) \cite{Labarthe2014, Roney2019}.
The transmural fibre orientation was set to be the endocardial and the epicardial orientation for elements below and above 50\% of the wall thickness, respectively.
We refer to \cite{Strocchi2023} for further details about this patient-specific geometry.

\section{Mathematical and numerical modeling of the 3D-0D solver}
\label{app:methods:3D0D}

Let $\Omega \subset \mathbb{R}^3$ be the domain corresponding to the patient-specific four-chamber heart.
We employ the reaction-Eikonal model without diffusion for cardiac electrophysiology \cite{Neic2017}.
We report the Eikonal model in Equation~\eqref{eq:eikonal}.
Given $\mathbf{V}(\mathbf{x})$ containing the squared local conduction velocities (CV) in the fibres, sheet and normal to sheet directions, and sites of initial activation $\Gamma$, this equation allows to find the local activation times $t_a(\mathbf{x})$ at node location $\mathbf{x}$, with initial activation occurring at a prescribed time $t_0$:

\begin{equation} \label{eq:eikonal}
    \left\{
    \begin{aligned}
        \sqrt{\nabla t_a(\mathbf{x})^{T} \mathbf{V}(\mathbf{x}) \nabla t_a(\mathbf{x})} &= 1 \qquad \mathbf{x} \in \Omega, \\
        t_a(\mathbf{x}) &= t_0 \qquad \mathbf{x} \in \Gamma. \\
    \end{aligned}
    \right.
\end{equation}

\begin{figure}[t!]
    \centering
    \includegraphics[width=\textwidth]{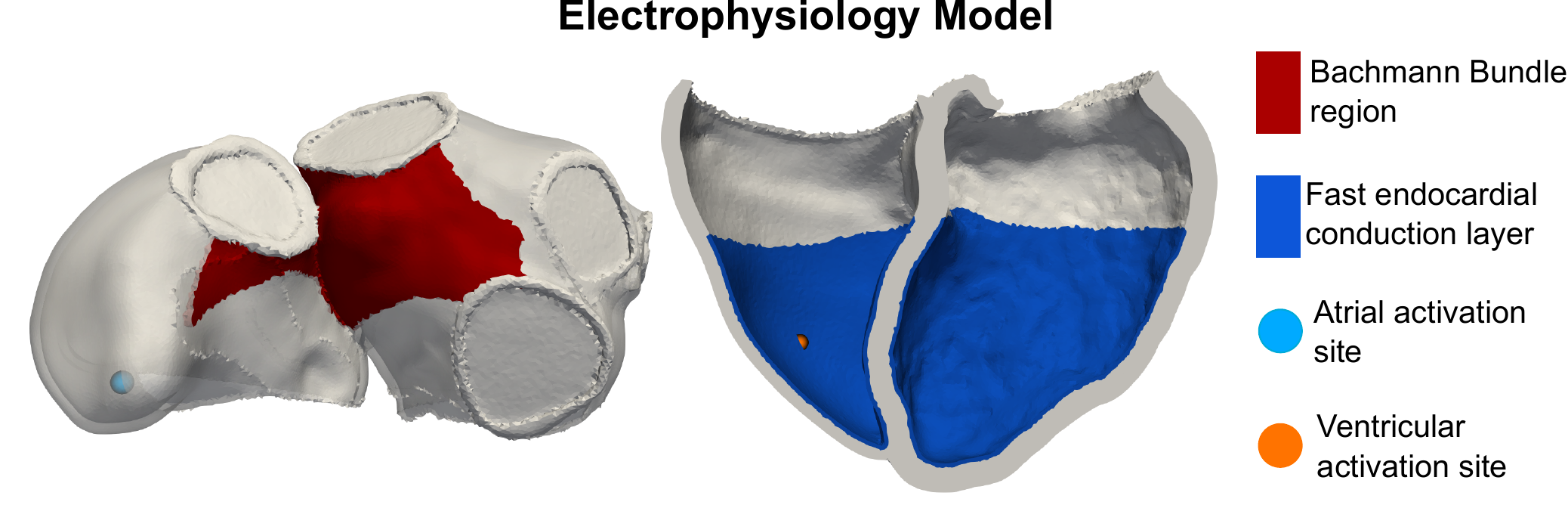}
    \caption{Electrophysiology simulations. Atria with the region representing the Bachmann bundle and the atrial activation site (left). Ventricles with the fast endocardial conduction layer and the ventricular activation site (right) \cite{Strocchi2023}.}
    \label{fig:EPmodel}
\end{figure}

We represent atria and ventricles as transversely isotropic conductive regions.
In particular, we assign CVs in the fibre direction (CV\textsubscript{f,V} and CV\textsubscript{f,A}) and anisotropy ratios (k\textsubscript{ft,V} and k\textsubscript{ft,A}), respectively.
The remaining regions are considered as passive.
To represent fast endocardial activation due to the His--Purkinje system, we introduce a 1-mm element thick endocardial layer extending up to 70\% in the apico-basal direction of the ventricles \cite{Lee2019, Strocchi2020Cohort}, with faster CV compared to the rest of ventricular myocardium of a factor $k_{\text{FEC}}$ (Figure~\ref{fig:EPmodel}, right).
We account for the Bachmann bundle by defining a region between the left atrium (LA) and the right atrium (RA) with fast CV compared to the rest of the atrial myocardium of a factor $k_{\text{BB}}$ (Figure~\ref{fig:EPmodel}, left) \cite{Roney2019}.
To fully control the atrioventricular (AV) delay, we define a passive region along the AV plane to insulate the atria from the ventricles.
Atrial activation is triggered at the location of the RA lead, while ventricular activation is initiated at the RV lead location with a delay defined by the AV delay, included as a free parameter in the simulator (AV\textsubscript{delay}).
The RA and RV lead locations were selected by segmenting the pacemaker leads from the CT image by thresholding the image intensity.

We employ the Courtemanche-Ramirez-Nattel (CRN) \cite{Courtemanche1998} and the Tomek-Rodriguez-O'Hara-Rudy (ToR-ORd) ionic model with dynamic intracellular chloride \cite{Tomek2019} for atrial and ventricular cardiomyocytes, respectively.
We induce the initial increase in the transmembrane potential $V_{\text{m}}$ by imposing a foot current that acts as a local stimulus, activating the cell membrane in each point $\mathbf{x}$ of the domain $\Omega$ at the local activation time $t_a(\mathbf{x})$ computed with the Eikonal model \cite{Neic2017}.

The intracellular calcium transient obtained from the ionic model is provided as an input to the Land contraction model \cite{Land2017} to compute the active tension transient in atria and ventricles.
For the sake of simplicity, we assume that active contraction occurs in the fibre direction only.
Prior to the 3D-0D closed-loop electromechanical simulations, the ToR-ORd-Land and CRN-Land cell models were run for 500 heartbeats at a basic cycle length of $0.854 \; \si{\second}$, which corresponds to the heartbeat period $\THB$ of the patient, to reach a steady state.

We use the transversely isotropic Guccione model for atrial and ventricular passive mechanics \cite{Guccione1991_2}, according to which the strain energy function takes the following expression:

\begin{align} \label{eq:Guccione}
    \Psi(\mathbf{E}) &= \frac{a}{2}\left[e^Q-1\right]+\frac{\kappa}{2}\left(\log J\right)^2 \\
    Q &= b_{\text{f}}E_{\text{ff}}^2 + 2b_{\text{ft}}(E_{\text{fs}}^2+E_{\text{fn}}^2)+b_{\text{t}}(E_{\text{ss}}^2+E_{\text{nn}}^2+2E_{\text{sn}}^2)\;,\notag
\end{align}

\noindent where $J$ is the determinant of the deformation tensor, $\mathbf{E}$ represents the Green-Lagrange strain tensor and $f$, $s$ and $n$ are the fibre, sheet and normal to sheet directions.
$a$, $b_{\text{f}}$, $b_{\text{ft}}$ and $b_{\text{t}}$ are the stiffness parameters, whereas $\kappa = 1000$ kPa is the bulk modulus, penalising volume changes and therefore enforcing quasi-incompressibility \cite{Flory1961, Ogden1978}.
Passive material properties of all the other cardiac tissues are represented by means of a Neo-Hookean model, with the stiffness parameters following previous studies \cite{Strocchi2020Cohort, Strocchi2020Simulating}. 

\begin{figure}[t!]
    \centering
    \includegraphics[width=\textwidth]{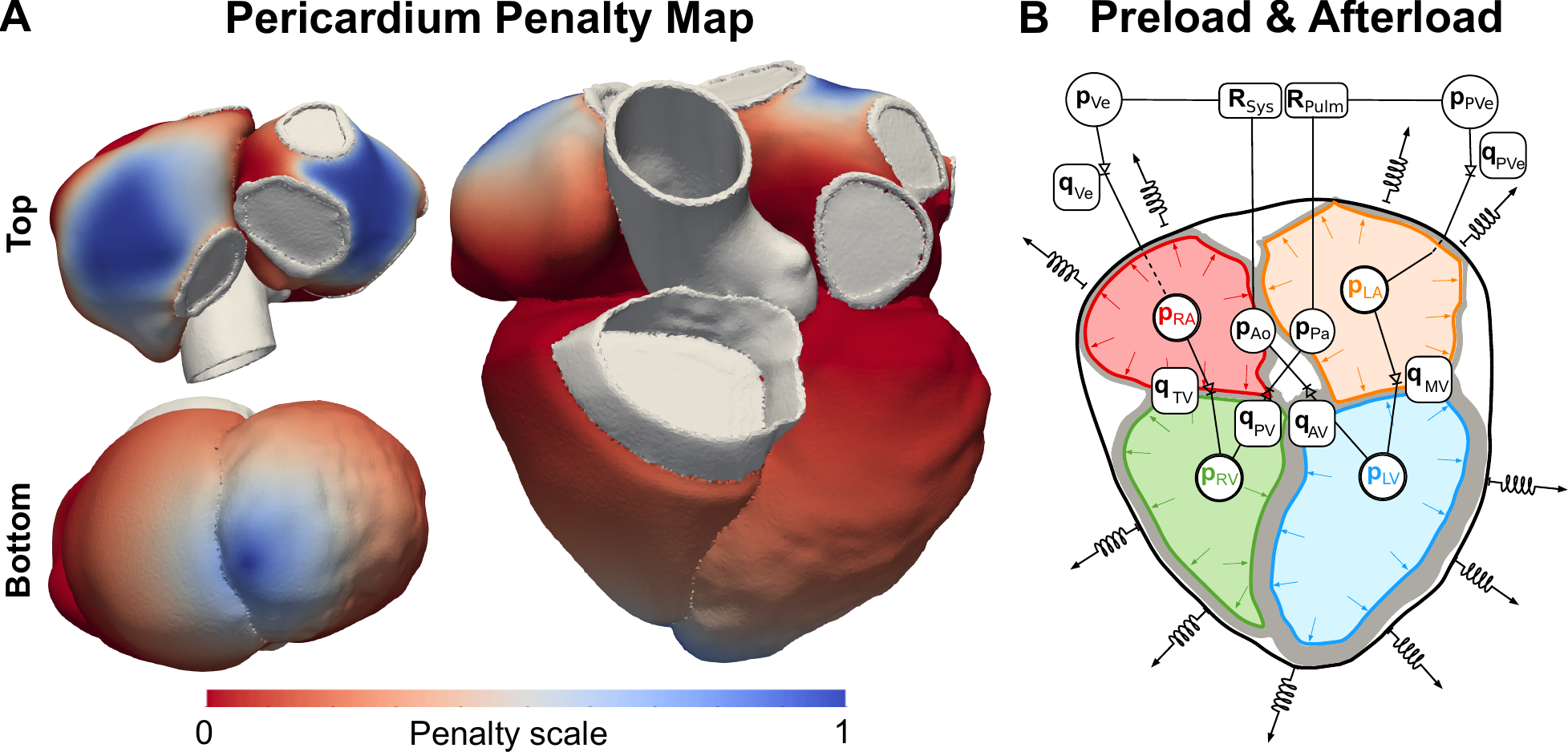}
    \caption{Boundary conditions for the mechanical model. \textbf{A} Penalty map scaling the normal spring stiffness for the effect of the pericardium. \textbf{B} Afterload and preload boundary conditions represented with CircAdapt. Symbols and abbreviations: \textbf{p}=pressure, \textbf{R}=resistance, \textbf{q}=flow across a valve, LV=left ventricle, RV=right ventricle=, LA=left atrium, RA=right atrium, Ao=aorta, Pa=pulmonary artery, Ve=veins, PVe=pulmonary veins, sys=systemic, pulm=pulmonary, MV=mitral valve, TV=tricuspid valve, AV=aortic valve, PV=pulmonary valve \cite{Strocchi2023}.}
    \label{fig:BCs}
\end{figure}

As described in \cite{Strocchi2020Cohort, Augustin2021}, we simulate the pericardium effect on the heart with normal springs with stiffness $k_{\text{peri}}$.
This value is scaled on the ventricles according to a map derived from motion data \cite{Strocchi2020Simulating}, to constrain the motion of the apex but not the base, allowing for physiological AV plane downward displacement during ventricular systole.
A similar analysis on the atria, described in \cite{Strocchi2021}, showed that the roof of the atria moved the least, while the regions around the AV plane moved the most, as they are stretched down by the contracting ventricles.
We therefore define a scaling map on the atria to include this constraint in the model, by assigning maximum penalty to the roof of the atria and zero penalty towards the AV plane (Figure~\ref{fig:BCs}A).
In addition, we apply omni-directional springs to the right inferior and superior pulmonary veins and at the superior vena cava rings.
The stiffness of these springs is fixed to 1.0 kPa/$\mu$m \cite{Strocchi2023}. 

The 3D four-chamber electromechanical model is coupled with the 0D closed-loop CircAdapt model \cite{Arts2005,Walmsley2015} (Figure~\ref{fig:BCs}B), which represents the following components of the circulatory system: aorta, pulmonary artery, veins, systemic and pulmonary peripheral resistances, the four cardiac valves (aortic, pulmonary, mitral and tricuspid) and flows across the pulmonary veins into the LA and across the systemic veins into the RA.
The monolithic 3D-0D coupling method is described in \cite{Augustin2021}.
Briefly, the pressures of the LA, LV, RA and RV were included as additional unknowns to the monolithic scheme, and the following equations are added to the equations of passive mechanics:

\begin{align}
\notag
    V^{\text{3D}}_{\text{LV}}(\mathbf{u},t) - V^{\text{0D}}(p_{\text{LV}},t) &= 0\;, \\ \notag
    V^{\text{3D}}_{\text{RV}}(\mathbf{u},t) - V^{\text{0D}}(p_{\text{RV}},t) &= 0\;, \\ \notag
    V^{\text{3D}}_{\text{LA}}(\mathbf{u},t) - V^{\text{0D}}(p_{\text{LA}},t) &= 0\;, \\ \notag
    V^{\text{3D}}_{\text{RA}}(\mathbf{u},t) - V^{\text{0D}}(p_{\text{RA}},t) &= 0\;,    \notag
\end{align}

\noindent where $V^{\text{3D}}$ and $V^{\text{0D}}$ are the volumes of the cavity computed from the deforming 3D mesh and predicted by the 0D model, respectively, $t$ is the time and $\mathbf{u}$ is the displacement field.

The ventricles of the end-diastolic mesh are unloaded from an end-diastolic left ventricle (LV) and right ventricle (RV) pressure, while the atria are not unloaded, under the assumption that the active tension in the atrial myocardium balances the pressure \cite{Land2017}.
During the unloading phase, we do not apply pericardial boundary conditions at the epicardium.
Then, prior to the start of the 3D-0D coupled simulation, we reloaded the ventricles to retrieve the end-diastolic mesh while the atrial pressure is initialised at 0 mmHg.
The electromechanical simulations always start at end-diastole and the pericardial boundary conditions are actived in this phase.
We remark that, differently from \cite{Strocchi2023}, the end-diastolic pressures do not act as additional model parameters but are instead prescribed as initial conditions for the system of LNODEs.

\begin{figure}[t!]
    \centering
    \includegraphics[width=0.75\textwidth]{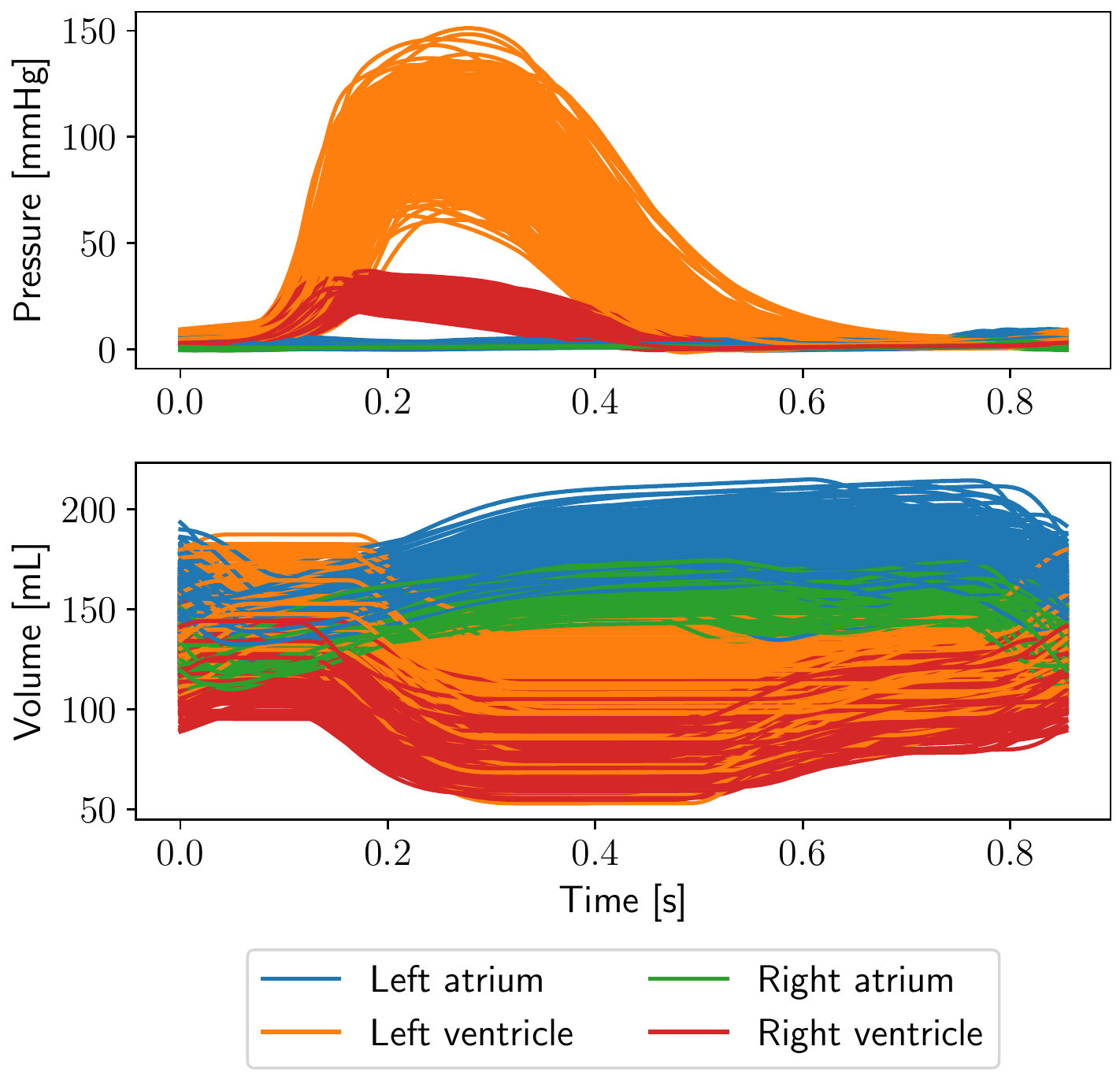}
    \caption{Pressure and volume time traces of the $\numSims = 405$ electromechanical simulations considered for training, validation and testing phases. We run the numerical simulations for 5 heartbeats (heartbeat period $\THB = 0.854 \; \si{\second}$) and then we take the last cardiac cycle.}
    \label{fig:pressuresvolumes}
\end{figure}

\begin{figure}[t!]
    \centering
    \includegraphics[width=0.75\textwidth]{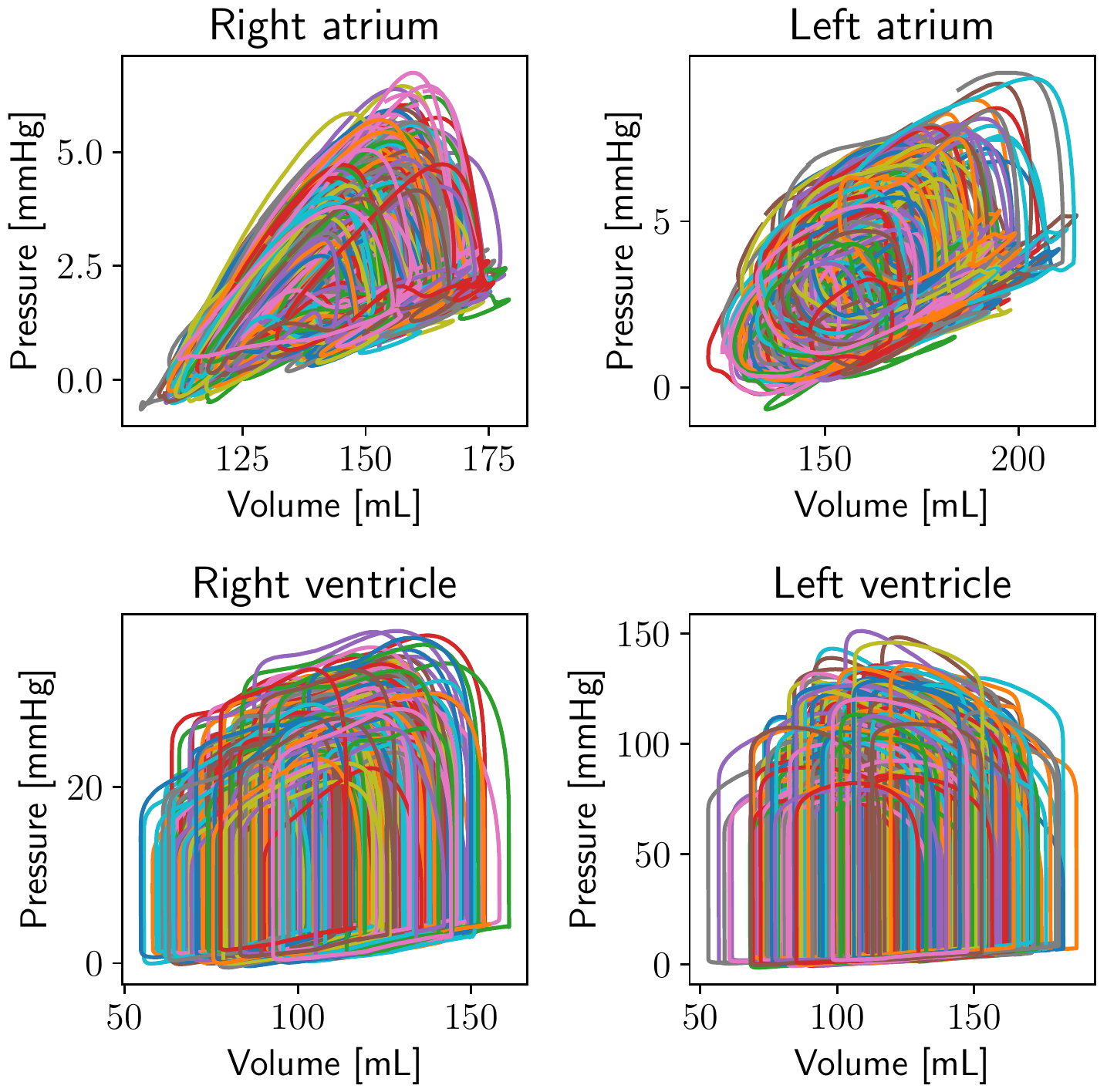}
    \caption{Corresponding pressure-volume loops of the $\numSims = 405$ electromechanical simulations used for training, validation and testing phases.}
    \label{fig:pvloops}
\end{figure}

To minimise the effect of these initial conditions, we run all numerical simulations for 5 heartbeats, to reach a near-to-steady-state behaviour, on a supercomputer endowed with 512 cores.
Figures~\ref{fig:pressuresvolumes} and~\ref{fig:pvloops} show the pressure-volume dynamics of all the $\numSims = 405$ electromechanical simulations considered for training, validation and testing phases.
Given the significant amount of required computational power, we select the linear and non-linear solver relative tolerances for passive mechanics in such a way to reduce the overall computational time while preserving accuracy \cite{Strocchi2023}.
In particular we set the maximum number of Newton iterations to 1 for the first three heartbeats. Indeed, as shown in \cite{Augustin2021}, this approach brings the numerical simulation closer to a steady state before solving nonlinear passive mechanics more accurately with more Newton iterations.
As a matter of fact, we set the maximum number of Newton iterations to 2 for the last two heartbeats in order to have a better approximation of the stretch rate for the cell model.
We also increase the tolerance for the numerical solution of the linearised system to 10\textsuperscript{-4}, for all heartbeats.
We show in \cite{Strocchi2023} that these numerical settings have limited effects on the pressure-volume dynamics simulated by the 3D-0D closed-loop electromechanical model while allowing for a 3 times speed-up in the total computational time.
We refer to \cite{Strocchi2023} for further details about the mathematical and numerical model.

\section{Model parameters}
\label{app:methods:parameters}
In Table~\ref{tab:parameterspace}, we report the list of parameters covering the whole cardiac and cardiovascular function that has been used to train the system of LNODEs.
The choice of the specific parameter values and their ranges is driven by the comprehensive study performed in \cite{Strocchi2023}, where Strocchi et al. train several Gaussian Processes Emulators to carry out Global Sensitivity Analysis and History Matching.
In particular, starting from 117 model parameters of interest, this technique allows to exclude unimportant ones, whereas the latter permits to identify implausible areas that would provide unphysiological outputs.

\begin{table}[p!]
    \begin{center}
        \hspace*{-2.5cm}
        \begin{tabular}{ l l l l}
            \toprule
            Parameter & Description & Range & Refs. \\
            \midrule
            $\IupmaxCRN$            & AT maximum $Ca^{2+}$ uptake rate into the sarcoplasmic reticulum network & [0.0028, 0.0080]     & \cite{Courtemanche1998} \\
            $\TrpnmaxCRN$           & AT total troponin C concentration in cytoplasm                           & [0.0359, 0.1039]     & \cite{Courtemanche1998} \\
            $\gCaLCRN$              & AT conductance of L-type $Ca^{2+}$ current                               & [0.08910, 0.1979]    & \cite{Courtemanche1998} \\
            $\TrefCRNLand$          & AT reference isometric tension                                           & [80.1142, 119.934]   & \cite{Land2017} \\
            $\permfiftyCRNLand$     & AT calcium/troponin complex when $50\%$ of crossbridges are blocked      & [0.1802, 0.5233]     & \cite{Land2017} \\
            $\npermCRNLand$         & AT Hill coefficient for $Ca^{2+}$-troponin and unbound sites             & [2.54302, 7.3862]    & \cite{Land2017} \\
            $\TRPNnCRNLand$         & AT $Ca^{2+}$-troponin cooperativity                                      & [1.0222, 2.9841]     & \cite{Land2017} \\
            $\drCRNLand$            & AT steady-state duty ratio                                               & [0.1271, 0.3711]     & \cite{Land2017} \\
            $\wfracCRNLand$         & AT steady-state ratio between weakly and strongly bound sites            & [0.2612, 0.7461]     & \cite{Land2017} \\
            $\TOTACRNLand$          & AT scale for distortion due to velocity of contraction                   & [12.7021, 37.1894]   & \cite{Land2017} \\
            $\phiCRNLand$           & AT distortion decay                                                      & [1.1738, 3.4305]     & \cite{Land2017} \\
            $\cafiftyCRNLand$       & AT reference $Ca^{2+}$ sensitivity                                       & [0.5678, 1.2880]     & \cite{Land2017} \\
            $\muCRNLand$            & AT scaling factor for weakly to strongly transition rate                 & [4.6065, 13.3936]    & \cite{Land2017} \\
            $\CVatria$              & Atrial conduction velocity in the fibre direction                        & [0.7508, 1.0269]     & \cite{Strocchi2023} \\
            $\kBB$                  & Bachmann bundle scaling factor                                           & [1.7011, 5.6372]     & \cite{Strocchi2023} \\
            $\aatria$               & AT bulk myocardium stiffness                                             & [1.5095, 2.4999]     & \cite{Land2018,Oken1957} \\
            $\bfatria$              & AT stiffness in the fibre direction                                      & [4.0493, 11.9966]    & \cite{Land2018,Nasopoulou2017} \\
            $\btatria$              & AT stiffness in the transverse plane                                     & [1.5192, 4.4989]     & \cite{Land2018,Nasopoulou2017} \\
            $\PCabToRORd$           & VE conductance of the background $Ca^{2+}$ current                       & [6.3660e-05, 1e-04]  & \cite{Tomek2019} \\
            $\trpnmaxToRORd$        & VE maximum troponin C concentration                                      & [0.065, 0.1228]      & \cite{Tomek2019} \\
            $\GncxbToRORd$          & VE conductance of the $Na^{2+}$-$Ca^{2+}$ exchanger                      & [0.0009, 0.0024]     & \cite{Tomek2019} \\
            $\TrefToRORdLand$       & VE reference isometric tension                                           & [127.846, 199.575]   & \cite{Land2017} \\
            $\permfiftyToRORdLand$  & VE calcium/troponin complex when $50\%$ of crossbridges are blocked      & [0.1764, 0.5117]     & \cite{Land2017} \\
            $\npermToRORdLand$      & VE Hill coefficient for $Ca^{2+}$-troponin and unbound sites             & [1.8542, 3.0441]     & \cite{Land2017} \\
            $\TRPNnToRORdLand$      & VE $Ca^{2+}$-troponin cooperativity                                      & [1.8390, 2.9980]     & \cite{Land2017} \\
            $\drToRORdLand$         & VE steady-state duty ratio                                               & [0.1263, 0.3627]     & \cite{Land2017} \\
            $\wfracToRORdLand$      & VE steady-state ratio between weakly and strongly bound sites            & [0.2884, 0.7462]     & \cite{Land2017} \\
            $\TOTAToRORdLand$       & VE scale for distortion due to velocity of contraction                   & [12.6888, 37.392]    & \cite{Land2017} \\
            $\ktmunblockToRORdLand$ & VE transition rate from blocked to unblocked binding site                & [0.0123, 0.0315]     & \cite{Land2017} \\
            $\cafiftyToRORdLand$    & VE reference $Ca^{2+}$ sensitivity                                       & [0.4071, 1.0490]     & \cite{Land2017} \\
            $\muToRORdLand$         & VE scaling factor for weakly to strongly transition                      & [1.5216, 4.4844]     & \cite{Land2017} \\
            $\CVventricles$         & VE conduction velocity in the fibre direction                            & [0.3832, 0.7967]     & \cite{Strocchi2023} \\
            $\kFEC$                 & Fast endocardial layer scaling factor                                    & [1.3250, 8.3687]     & \cite{Strocchi2023} \\
            $\aventricles$          & VE bulk myocardium stiffness                                             & [0.5006, 1.4998]     & \cite{Land2018} \\
            $\btventricles$         & VE stiffness in the transverse plane                                     & [1.5042, 4.49251]    & \cite{Land2018,Nasopoulou2017} \\
            $\alvrv$                & Scaling factor for $\aventricles$ in RV vs. LV                           & [1.0055, 1.9995]     & \cite{Oken1957} \\
            $\Treflvrv$             & Scaling factor for $\TrefToRORdLand$ in RV vs. LV                        & [0.5009, 0.9956]     & \cite{Oken1957} \\
            $\AVdelay$              & Atrioventricular delay                                                   & [0.1, 0.2]           & \cite{Hyde2016} \\
            $\kperi$                & Pericardial normal springs stiffness                                     & [0.0005, 0.0019]     & \cite{Strocchi2020Simulating} \\
            $\Rsys$                 & Systemic resistance scaling factor                                       & [1.0017, 3.9937]     & \cite{Augustin2021,Walmsley2015} \\
            $\Rpulm$                & Pulmonary resistance scaling factor                                      & [1.0020, 3.9980]     & \cite{Augustin2021,Walmsley2015} \\
            $\Aol$                  & Length of the aorta                                                      & [300.478, 498.745]   & \cite{Augustin2021,Walmsley2015} \\
            $\kArt$                 & Stiffness of the aorta                                                   & [6.0118, 9.9894]     & \cite{Augustin2021,Walmsley2015} \\
            \bottomrule
        \end{tabular}
        \caption{Parameter space explored by model $\modFOM$ and used for the ANN training. AT: atrial, VE: ventricular.}
        \label{tab:parameterspace}
    \end{center}
\end{table}

\section{Training of the Latent Neural Ordinary Differential Equations}
\label{app:methods:LNODE}

\begin{table}[h!]
    \hspace*{-0.75cm}
    \begin{tabular}{c ccccc c}
        \toprule
        \multirow{2}{*}{LNODE} & \multicolumn{5}{c}{Hyperparameters} & Trainable parameters \\
               & layers & neurons & num. states & loss integr. step [$\si{\second}$] & weights reg. & \# param. \\
        \midrule
        tuning & $\{1 \; ... \; 7\}$ & $\{5 \; ... \; 50\}$ & $\{8 \; ... \; 12\}$ & $[10^{-3}, 0.1]$ & $[10^{-4}, 1]$ & \\
        final  & 3                   & 13                   & 8                    & 0.0285           & 0.023          & 1'178 \\
        \bottomrule
    \end{tabular}
    \caption{Hyperparameters ranges and selected values for the final training stage of LNODEs.}
    \label{tab:hyperparameters}  
\end{table}

\noindent We perform hyperparameters tuning by employing $K$-fold ($K = 10$) cross validation over 400 electromechanical simulations.
An optimal set of hyperparameters is automatically found by running the Tree-structured Parzen Estimator (TPE) Bayesian algorithm \cite{Bergstra2011, Optuna2019} while monitoring the generalization error reported in the main text (Section~\ref{sec:methods:LNODE}, Equation~\ref{eqn:loss}) during $K$-fold cross validation.
We early stop bad hyperparameters configurations by means of the Asynchronous Successive Halving (ASHA) scheduler \cite{Asha2020, Hyperband2017}.
We rely on the Ray Python distributed framework for the implementation of this hyperparameters tuner \cite{Ray2018}.
Different ANNs associated to different hyperparameters settings are simultaneously trained with Message Passing Interface (MPI) on 40 cores of a high-performance computing facility at MOX, Dipartimento di Matematica, Politecnico di Milano. 
We also exploit Hyper-Threading via Open Multi-Processing (OpenMP) to speed-up tensorial operations in Tensorflow \cite{Tensorflow2015}.
We consider an hypercube as a search space for the following hyperparameters: number of layers and neurons of the ANN, number of states $\NumANNState$, loss function $\mathcal{L}(\ANNState(t), \ANNStateTilde(t); \ANNparamTrained)$ integration step $\Delta t_\text{ref}$ and $L^2$ weights regularization $\iota$.
For each configuration of hyperparameters, we perform 1'000 iterations with the first-order Adam optimizer \cite{Kingma2014}, starting with a learning rate of $10^{-2}$, and then we continue the training stage with 10'000 iterations of the second-order BFGS optimizer \cite{Goodfellow2016}.
In this way, we exploit the stochastic behavior of the Adam optimizer to explore the landscape of local minima, and then we properly reach convergence by means of the BFGS optimizer.
The ANN is always initialized with a new set of weights provided by a Glorot uniform distribution and zero values for biases.
In Table~\ref{tab:hyperparameters}, we report the initial hyperparameters ranges for tuning and the final optimized values.

\section{Global sensitivity analysis}
\label{app:methods:GSA}

\begin{figure}
       \centering 
       \includegraphics[width=1.1\textwidth]{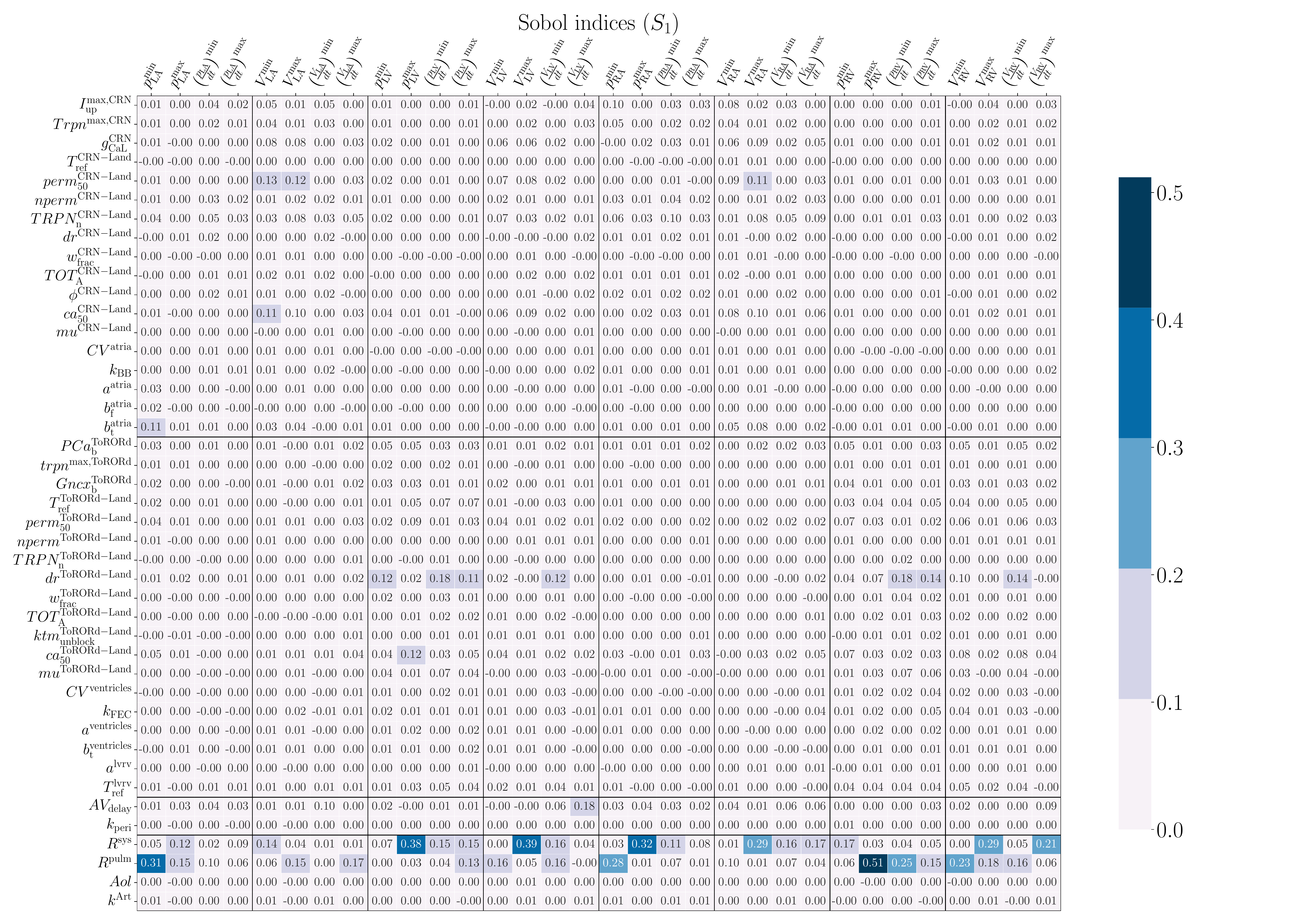}
       \caption{First-order Sobol indices $\SobolFirst{i}{j}$ computed by exploiting model $\modROM$.
                Each row corresponds to a parameter $\param$ of the 3D-0D closed-loop electromechanical model (see Table~\ref{tab:parameterspace}).
                Each column corresponds to a QoI $\qoi$ (maximum and minimum of the temporal traces of pressures or volumes and of their time derivatives).
                Both parameters and QoIs are split into a number of groups, separated by a black solid line.
                Specifically, from left to right, we list QoIs referred to LA, LV, RA and RV.
                From top to bottom, we list model parameters associated with atria, ventricles, whole-heart dynamics and cardiovascular system.}
       \label{fig:S1}
\end{figure}

To assess how much each model parameter $\param_i$ affects a pressure-volume biomarker of clinical interest for the atrial or ventricular function, that is a QoI $\qoi_j$, we compute Sobol indices via a probabilistically-driven variance-based sensitivity analysis \cite{sobol1990sensitivity}.
The first-order Sobol index $\SobolFirst{i}{j}$ evaluates the impact that a single parameter $\param_i$ has on a certain QoI $\qoi_j$, whereas the total-effect Sobol index $\SobolTotal{i}{j}$ also accounts for the interactions among parameters:

\begin{equation*}
    \SobolFirst{i}{j} = \frac{\variance_{\param_i}\left[ \expected_{\param_{\sim i}}\left[ \qoi_j | \param_i\right] \right]}{\variance\left[ \qoi_j \right]},
\;\;\;\;\;\;\;
    \SobolTotal{i}{j} = \frac{\expected_{\param_{\sim i}}\left[ \variance_{\param_i}\left[ \qoi_j | \param_{\sim i}\right] \right]}{\variance\left[ \qoi_j \right]}
                  = 1 - \frac{\variance_{\param_{\sim i}}\left[ \expected_{\param_i}\left[ \qoi_j | \param_{\sim i}\right] \right]}{\variance\left[ \qoi_j \right]},
\end{equation*}
where $\param_{\sim i}$ indicates the set of all parameters excluding the $i^{th}$ one.

We employ the Saltelli's method and model $\modROM$ to estimate these Sobol indices \cite{homma1996importance,saltelli2002making}.
This allows for a linear increase in the number of samples $\NumSamples$ with respect to the number of parameters $\NumParams$ once a certain accuracy is prescribed.
Specifically, the number of samples $\NumSamples$ in the parameter space scales as $N (2 \NumParams + 2)$, being $N$ a user defined value.
In this work, we set $N = 8000$, for a total of 704'000 samples, that allows for small confidence intervals around the first-order and total-effect Sobol indices, respectively.

The model parameters may arbitrarily vary in the training ranges defined in Table~\ref{tab:parameterspace}.
The QoIs are given by the maximum and minimum values of the four-chamber pressures, volumes and corresponding time derivatives of the simulated heartbeat with the trained LNODEs.
We employ the forward Euler method with a fixed time step $\Delta t = 10^{-3} \; \si{\second}$ for all the numerical simulations with model $\modROM$.

We report first-order Sobol indices in Figure~\ref{fig:S1}.
By comparing them to total-effect Sobol indices in Figure~\ref{fig:ST} from the main text, we notice that there are some significant differences.
This implies that the variation of the individual model parameters is less dominant than their high-order interactions.
This holds for most of the influential parameters coming from the CRN-Land model, such as $\gCaLCRN$, $\permfiftyCRNLand$, $\TRPNnCRNLand$ and $\cafiftyCRNLand$, as well as driven from the ToRORd-Land model, such as $\permfiftyToRORdLand$, $\TrefToRORdLand$ and $\cafiftyToRORdLand$. 
Similar considerations can be done for other ventricular ($\CVventricles$, $\kFEC$) and cardiovascular parameters ($\Rsys$, $\Rpulm$) with respect to specific QoIs.
As in the case with total-effect Sobol indices, we notice that the QoIs associated with a given area of the cardiovascular system are still mostly determined by the parameters associated with the same region.
However, similar important exceptions can be outlined, which is the case for the systemic and pulmonary resistances ($\Rsys$, $\Rpulm$), which explain variability for the whole-heart.
The AV delay $\AVdelay$ has a major impact on almost all biomarkers and can be inferred from the patient-specific volume traces over time without the need for parameter calibration.

Finally, we remark that Sobol indices are affected by the amplitude of the ranges in which the parameters are varied.
In particular, the wider the range associated with a parameter, the greater the associated Sobol indices will be, as the parameter in question potentially generates greater variability in the QoI.
Therefore, we stress that the results shown here and in the main text are valid for the specific ranges we used.

\section{Robust parameter estimation}
\label{app:methods:PE}

Let $\param \in \paramSpace \subset \mathbb{R}^{\NumParams}$ be a subset of parameters for the model $\modFOM$ that we want to calibrate with model $\modROM$, when certain time-dependent QoIs are provided as observations.
We carry out the Maximum a Posteriori (MAP) estimation by solving a constrained optimal control problem \cite{Salvador2023}:
\begin{equation}
\min_{\param \in \paramSpace} J(\param),
\label{eqn: OCP}
\end{equation}
with the following cost functional:
\begin{equation}
\begin{split}
J(\param) &= \sum_{\text{i} \in \{ \text{LA, LV, RA, RV} \}} \epsilon_\text{i} \dfrac{||p_\text{i}(t;\param) - \hat{p}_\text{i}(t)||_{\text{L}^2}^2}{\mu_{p_\text{i}^2}} \\
& + \sum_{\text{i} \in \{ \text{LA, LV, RA, RV} \}} \zeta_\text{i} \dfrac{||V_\text{i}(t;\param) - \hat{V}_\text{i}(t)||_{\text{L}^2}^2}{\mu_{V_\text{i}^2}} \\
\end{split}
\label{eqn: costfunctional}
\end{equation}
where $t \in [T - \THB, T]$, that is by considering the last heartbeat.

Pressures $p_\text{i}$ and volumes $V_\text{i}$ are numerical solutions of model $\modROM$, while observations $\hat{p}_\text{i}$ and $\hat{V}_\text{i}$ may come from either in silico numerical simulations or clinical data.
In particular, in this work we only focus on observations coming from model $\modFOM$.
Coefficients $\epsilon_\text{i}$ to $\zeta_\text{i}$ weigh the pressure and volume traces over time.
In this paper, we balance the different terms in Equation~\eqref{eqn: costfunctional}, i.e. we set $\epsilon_\text{i}$ and $\zeta_\text{i}$ to be either 0 or 1 according to the specific test case.
The normalization terms $\mu_{p_\text{i}^2}$ and $\mu_{V_\text{i}^2}$ are defined by averaging the squared values of pressure and volume traces over time, for $t \in [T - \THB, T]$.

For the sake of clarity, we recall the mathematical formulation for a system of LNODEs here below:
\begin{equation}
    \left\{
    \begin{aligned}
        \frac{d \ANNState(t)}{d t} &= \ANNRhs\left(
            \ANNState(t),
            \cos\left(\tfrac{2 \pi (t - \AVdelay)}{\THB}\right),
            \sin\left(\tfrac{2 \pi (t - \AVdelay)}{\THB}\right),
            \param;
            \ANNparam
        \right)
        && \text{for } t \in (0, \THB],\\
        \ANNState(0) &= \ANNState_{0}.  && \\
    \end{aligned}
    \right.
    \label{eqn: forwardproblem}
\end{equation}
During each iteration of the optimization problem, we compute $\tfrac{d J}{d \param}$ to minimize the loss function $J = J(\param)$.
We solve Equation~\eqref{eqn: forwardproblem} forward in time, i.e. for $t \in (T - \THB, T]$, using an ODE numerical solver.
Then, we solve an adjoint ODE system backward in time by exploiting reverse-mode differentiation \cite{Chen2019}:
\begin{equation} 
\dfrac{d \adjoint(t)}{d t} = -\adjoint(t)^T \dfrac{\partial \mathbf{\ANNRhs}}{\partial \ANNState} \left(\ANNState(t), \cos\left(\tfrac{2 \pi (t - \AVdelay)}{\THB}\right), \sin\left(\tfrac{2 \pi (t - \AVdelay)}{\THB}\right), \param; \ANNparam \right),
\label{eqn: adjointproblem}
\end{equation}
being $\adjoint(t) = dJ/d\ANNState(t)$ the adjoint state.

Finally, the gradient of $J$ with respect to $\param$ reads \cite{Chen2019}:
\begin{equation}
\frac{d J}{d \param} = - \int_T^0 \adjoint(t)^T \frac{\partial \mathbf{\ANNRhs}}{\partial \param}\left(\ANNState(t), \cos\left(\tfrac{2 \pi (t - \AVdelay)}{\THB}\right), \sin\left(\tfrac{2 \pi (t - \AVdelay)}{\THB}\right), \param; \ANNparam \right) dt
\label{eqn: lossgradient}
\end{equation}
All vector-jacobian products in Equations~\ref{eqn: adjointproblem} and~\ref{eqn: lossgradient} are evaluated using matrix-free methods, automatic differentiation and automatic vectorization.
These overall define an efficient numerical strategy accounting for very small memory requirements \cite{Africa2023}.

We use the Limited-memory Broyden-Fletcher-Goldfarb-Shanno (L-BFGS) algorithm to solve the optimal control problem~\eqref{eqn: OCP} \cite{Liu1989}.
We employ the forward Euler method with a fixed time step $\Delta t = 10^{-3} \; \si{\second}$ to solve Equations~\ref{eqn: forwardproblem} and~\ref{eqn: adjointproblem} at each L-BFGS iteration \cite{Dormand1980}.
We remark that the optimization process is constrained according to the model parameters ranges reported in Table~\ref{tab:parameterspace}.


Once we provide pointwise values of model parameters $\param_\text{MAP}$ via MAP estimation, we evaluate the uncertainty of these estimated values by means of Hamiltonian Monte Carlo (HMC) \cite{Betancourt2017}.
This method for inverse uncertainty quantification (UQ) allows to find an approximation of either the marginal or joint posterior distribution $\mathbb{P}(\param | \mathbf{x})$ over $\param$.

Let $\momentum \in \mathbb{R}^{\NumMomentum}$ be a vector containing auxiliary momentum variables.
We define the conditional probability distribution of $\momentum$ given $\param$ as \cite{Betancourt2017}:
\begin{equation*}
\mathbb{P}(\momentum, \param) = \mathbb{P}(\momentum | \param) \mathbb{P}(\param),
\end{equation*}
being $\mathbb{P}(\param)$ the prior probability distribution over $\param$.
Then, by employing the kinetic energy $\mathbb{K}(\momentum | \param)=-\log \mathbb{P}(\momentum | \param)$ and the potential energy $\mathbb{U}(\param)=-\log \mathbb{P}(\param)$, we introduce the Hamiltonian function \cite{Betancourt2017}:
\begin{equation*}
\mathbb{H}(\momentum, \param) = -\log \mathbb{P}(\momentum, \param) = -\log \mathbb{P}(\momentum | \param) -\log \mathbb{P}(\param) = \mathbb{K}(\momentum | \param) + \mathbb{U}(\param).
\end{equation*}
Finally, we solve a coupled system of ODEs in $(\param, \momentum)$ to  advance the value of the parameters vector $\param=\param(t)$ from its current state:
\begin{equation} 
\begin{cases}
\dfrac{d \param}{d t} = \dfrac{\partial \mathbb{H}}{\partial \momentum} & \qquad \text{for } t \in (0, \overline{T}], \\
\dfrac{d \momentum}{d t} = -\dfrac{\partial \mathbb{H}}{\partial \param} & \qquad \text{for } t \in (0, \overline{T}],
\end{cases}
\label{eqn: HMC}
\end{equation}
where $t$ represents a fictitious time variable in the parametric space for $\param$.
We solve Equation~\ref{eqn: HMC} by means of the leapfrog time scheme and we employ the No-U-Turn Sampler (NUTS) extension of HMC, so that the number of virtual time steps is automatically determined and not user-defined \cite{Homan2014}.

We fix $\overline{\Delta t}=10^{-3}$ and we perform 750 iterations, for all the test cases.
Among them, the first 250 iterations consist of an initial burn-in phase and are not retained for the approximation of the posterior distribution $\mathbb{P}(\param | \mathbf{x})$.
We initialize the NUTS sampler by considering $\mathbb{P}(\param) \sim\ U(\param_\text{MAP} - \chi \param_\text{MAP}, \param_\text{MAP} + \chi \param_\text{MAP})$, being $\chi = 0.1$ a suitable parameter to define a uniform prior distribution around the MAP estimation $\param_\text{MAP}$.
With respect to standard Markov Chain Monte Carlo (MCMC), where multiple chains are usually required to achieve proper convergence to the posterior distribution \cite{Regazzoni2022EMROM}, here we will always run a single chain, as this is sufficient to provide meaningful results.
This is also motivated by the suitable initialization of the parameters, which is related to the MAP estimation $\param_\text{MAP}$.
We declare convergence when the Gelman-Rubin diagnostic provides a value less than 1.1 for all model parameters and there are no divergent transitions \cite{Vats2018}.

We account for the surrogate modeling error during the parameter identification process, as all the test cases of this paper are based on time-dependent QoIs coming from model $\modFOM$.
We consider normal distributions around the estimated values of these QoIs, i.e. $\mathcal{N}(p_\text{i}(t;\param), \NoiseCov_{\text{ANN,i}})$ and $\mathcal{N}(V_\text{i}(t;\param), \NoiseCov_{\text{ANN,i}})$, for $\text{i} \in \{ \text{LA, LV, RA, RV} \}$.
We introduce a zero-mean Gaussian process $\mathcal{G}\mathcal{P}(\boldsymbol{0}, k(\boldsymbol{t}, \boldsymbol{t}'))$, where $k(\boldsymbol{t}, \boldsymbol{t}')=\sigma^2 \exp \left( \tfrac{-||\boldsymbol{t}-\boldsymbol{t}'||^2}{2 \lambda^2} \right)$ is the exponentiated quadratic kernel \cite{Rasmussen2005}.
Amplitude $\sigma$ is independently computed for all the relevant pressure and volume time traces by looking at the pointwise differences of the outputs computed with model $\modROM$ and model $\modFOM$ \cite{Regazzoni2022EMROM}.
This leads to $\sigma_{\text{p}_\text{LA}} = 0.13 \; \si{\mmHg}$, $\sigma_{\text{p}_\text{LV}} = 2.30 \; \si{\mmHg}$, $\sigma_{\text{p}_\text{RA}} = 0.09 \; \si{\mmHg}$, $\sigma_{\text{p}_\text{RV}} = 0.46 \; \si{\mmHg}$, $\sigma_{\text{V}_\text{LA}} = 1.82 \; \si{\milli\liter}$, $\sigma_{\text{V}_\text{LV}} = 1.34 \; \si{\milli\liter}$, $\sigma_{\text{V}_\text{RA}} = 2.29 \; \si{\milli\liter}$ and $\sigma_{\text{V}_\text{RV}} = 1.38 \; \si{\milli\liter}$.
The correlation length $\lambda$ is estimated with 1'000 Adam iterations that minimize the negative log likelihood of the observed surrogate modeling error \cite{Kingma2014}.
We consider a unique value of $\lambda = 0.02$, as we observe similar correlation lengths for all the time-dependent QoIs.
The full covariance matrix $\NoiseCov_{\text{ANN,i}}$ can be then generated by means of the tuned kernel function:
\begin{equation*}
\NoiseCov_{\text{ANN,i}}(t_\text{j}, t_\text{k}) = \sigma_{\text{i}}^2 \exp \left[ \dfrac{-(t_\text{j} - t_\text{k})^2}{2 \lambda^2} \right] \;\;\; \text{for} \;\;\; \text{i} \in \{ \text{LA, LV, RA, RV} \},
\end{equation*}
being $t_\text{j}$ and $t_\text{k}$ discrete time points in $[T - \THB, T]$.
We remark that additive measurement errors driven by instrument sensitivities, surrounding environment and human intervention, related to noisy (realistic) observations, may be easily incorporated in our UQ framework as well \cite{Salvador2023}.

In Tables~\ref{tab:TLV}-\ref{tab:Tall}, we report the true values and mean plus/minus two standard deviations for all the estimated parameters in the different test cases ($\TLV$, $\Tventricles$, $\Tatria$, $\Tall$), for each numerical simulation of the testing set.
We notice that the true parameter value is always properly captured in the range of uncertainty of the corresponding estimation.

\begin{sidewaystable}[p!]
    \begin{center}
        \begin{tabular}{ l|ll|ll|ll|ll|ll }
            \toprule
            Parameter & \makecell{Ground \\ truth} & \makecell{Simulation \\ 1} & \makecell{Ground \\ truth} & \makecell{Simulation \\ 2} & \makecell{Ground \\ truth} & \makecell{Simulation \\ 3} & \makecell{Ground \\ truth} & \makecell{Simulation \\ 4} & \makecell{Ground \\ truth} & \makecell{Simulation \\ 5} \\
            \midrule
            $\drToRORdLand$        & 0.23 & 0.21 $\pm$ 0.06 & 0.17 & 0.17 $\pm$ 0.06 & 0.19 & 0.25 $\pm$ 0.06 & 0.13 & 0.13 $\pm$ 0.06 & 0.20 & 0.22 $\pm$ 0.06 \\ 
            $\cafiftyToRORdLand$   & 0.77 & 0.77 $\pm$ 0.18 & 0.53 & 0.66 $\pm$ 0.18 & 0.61 & 0.75 $\pm$ 0.18 & 0.57 & 0.61 $\pm$ 0.18 & 0.51 & 0.50 $\pm$ 0.18 \\
            $\Rsys$                & 3.28 & 3.28 $\pm$ 0.64 & 2.50 & 2.26 $\pm$ 0.63 & 3.84 & 3.21 $\pm$ 0.64 & 3.57 & 3.28 $\pm$ 0.63 & 2.26 & 2.77 $\pm$ 0.64 \\
            $\Rpulm$               & 2.63 & 2.74 $\pm$ 0.63 & 2.84 & 2.64 $\pm$ 0.63 & 3.12 & 2.80 $\pm$ 0.63 & 2.45 & 2.55 $\pm$ 0.63 & 3.40 & 3.02 $\pm$ 0.63 \\
            \bottomrule
        \end{tabular}
        \caption{$\TLV$: true value and mean plus/minus two standard deviations associated to the estimated model parameters for all the $\numTest = 5$ electromechanical simulations.}
        \label{tab:TLV}
    \end{center}

    \begin{center}
        \begin{tabular}{ l|ll|ll|ll|ll|ll }
            \toprule
            Parameter & \makecell{Ground \\ truth} & \makecell{Simulation \\ 1} & \makecell{Ground \\ truth} & \makecell{Simulation \\ 2} & \makecell{Ground \\ truth} & \makecell{Simulation \\ 3} & \makecell{Ground \\ truth} & \makecell{Simulation \\ 4} & \makecell{Ground \\ truth} & \makecell{Simulation \\ 5} \\
            \midrule
            $\drToRORdLand$        & 0.23 & 0.23 $\pm$ 0.05 & 0.17 & 0.16 $\pm$ 0.05 & 0.19 & 0.23 $\pm$ 0.05 & 0.13 & 0.12 $\pm$ 0.05 & 0.20 & 0.21 $\pm$ 0.05 \\ 
            $\cafiftyToRORdLand$   & 0.77 & 0.60 $\pm$ 0.26 & 0.53 & 0.53 $\pm$ 0.25 & 0.61 & 0.62 $\pm$ 0.26 & 0.57 & 0.46 $\pm$ 0.26 & 0.51 & 0.69 $\pm$ 0.26 \\
            $\permfiftyToRORdLand$ & 0.25 & 0.32 $\pm$ 0.14 & 0.26 & 0.36 $\pm$ 0.14 & 0.37 & 0.35 $\pm$ 0.14 & 0.27 & 0.35 $\pm$ 0.14 & 0.31 & 0.23 $\pm$ 0.14 \\
            $\Rsys$                & 3.28 & 3.30 $\pm$ 0.20 & 2.50 & 2.30 $\pm$ 0.21 & 3.84 & 3.73 $\pm$ 0.20 & 3.57 & 3.44 $\pm$ 0.20 & 2.26 & 2.38 $\pm$ 0.21 \\
            $\Rpulm$               & 2.63 & 2.67 $\pm$ 0.40 & 2.84 & 2.63 $\pm$ 0.41 & 3.12 & 2.82 $\pm$ 0.41 & 2.45 & 2.47 $\pm$ 0.41 & 3.40 & 3.15 $\pm$ 0.40 \\
            \bottomrule
        \end{tabular}
        \caption{$\Tventricles$: true value and mean plus/minus two standard deviations associated to the estimated model parameters for all the $\numTest = 5$ electromechanical simulations.}
        \label{tab:Tventricles}
    \end{center}
\end{sidewaystable}

\begin{sidewaystable}[p!]
    \begin{center}
        \begin{tabular}{ l|ll|ll|ll|ll|ll }
            \toprule
            Parameter & \makecell{Ground \\ truth} & \makecell{Simulation \\ 1} & \makecell{Ground \\ truth} & \makecell{Simulation \\ 2} & \makecell{Ground \\ truth} & \makecell{Simulation \\ 3} & \makecell{Ground \\ truth} & \makecell{Simulation \\ 4} & \makecell{Ground \\ truth} & \makecell{Simulation \\ 5} \\
            \midrule
            $\drToRORdLand$        & 0.23 & 0.20 $\pm$ 0.07 & 0.17 & 0.16 $\pm$ 0.07 & 0.19 & 0.18 $\pm$ 0.07 & 0.13 & 0.18 $\pm$ 0.07 & 0.20 & 0.21 $\pm$ 0.07 \\
            $\permfiftyToRORdLand$ & 0.25 & 0.29 $\pm$ 0.09 & 0.26 & 0.33 $\pm$ 0.09 & 0.37 & 0.45 $\pm$ 0.09 & 0.27 & 0.39 $\pm$ 0.09 & 0.31 & 0.34 $\pm$ 0.09 \\
            $\cafiftyCRNLand$      & 1.09 & 1.27 $\pm$ 0.23 & 0.53 & 1.16 $\pm$ 0.23 & 1.09 & 0.78 $\pm$ 0.23 & 1.06 & 0.92 $\pm$ 0.23 & 0.73 & 0.68 $\pm$ 0.23 \\
            $\TRPNnCRNLand$        & 2.89 & 2.14 $\pm$ 0.90 & 1.65 & 2.07 $\pm$ 0.90 & 2.14 & 2.37 $\pm$ 0.90 & 1.95 & 2.16 $\pm$ 0.90 & 2.08 & 2.08 $\pm$ 0.90 \\
            $\gCaLCRN$             & 0.13 & 0.11 $\pm$ 0.04 & 0.12 & 0.13 $\pm$ 0.04 & 0.19 & 0.16 $\pm$ 0.04 & 0.14 & 0.15 $\pm$ 0.04 & 0.13 & 0.17 $\pm$ 0.04 \\
            $\btatria$             & 3.19 & 3.96 $\pm$ 1.04 & 2.40 & 2.59 $\pm$ 1.04 & 2.24 & 1.77 $\pm$ 1.04 & 2.86 & 2.86 $\pm$ 1.04 & 2.50 & 1.86 $\pm$ 1.04 \\
            $\Rsys$                & 3.28 & 3.33 $\pm$ 0.38 & 2.50 & 2.28 $\pm$ 0.38 & 3.84 & 3.42 $\pm$ 0.38 & 3.57 & 3.61 $\pm$ 0.38 & 2.26 & 2.56 $\pm$ 0.38 \\
            $\Rpulm$               & 2.63 & 2.98 $\pm$ 0.67 & 2.84 & 2.62 $\pm$ 0.67 & 3.12 & 2.82 $\pm$ 0.67 & 2.45 & 1.92 $\pm$ 0.67 & 3.40 & 3.18 $\pm$ 0.67 \\
            \bottomrule
        \end{tabular}
        \caption{$\Tatria$: true value and mean plus/minus two standard deviations associated to the estimated model parameters for all the $\numTest = 5$ electromechanical simulations.}
        \label{tab:Tatria}
    \end{center}

    \begin{center}
        \begin{tabular}{ l|ll|ll|ll|ll|ll }
            \toprule
            Parameter & \makecell{Ground \\ truth} & \makecell{Simulation \\ 1} & \makecell{Ground \\ truth} & \makecell{Simulation \\ 2} & \makecell{Ground \\ truth} & \makecell{Simulation \\ 3} & \makecell{Ground \\ truth} & \makecell{Simulation \\ 4} & \makecell{Ground \\ truth} & \makecell{Simulation \\ 5} \\
            \midrule
            $\drToRORdLand$        & 0.23 & 0.27 $\pm$ 0.04 & 0.17 & 0.15 $\pm$ 0.04 & 0.19 & 0.21 $\pm$ 0.04 & 0.13 & 0.15 $\pm$ 0.04 & 0.20 & 0.22 $\pm$ 0.04 \\
            $\permfiftyToRORdLand$ & 0.25 & 0.23 $\pm$ 0.12 & 0.26 & 0.22 $\pm$ 0.12 & 0.37 & 0.30 $\pm$ 0.12 & 0.27 & 0.30 $\pm$ 0.12 & 0.31 & 0.18 $\pm$ 0.13 \\
            $\cafiftyToRORdLand$   & 0.77 & 0.85 $\pm$ 0.22 & 0.53 & 0.62 $\pm$ 0.22 & 0.61 & 0.70 $\pm$ 0.22 & 0.57 & 0.50 $\pm$ 0.22 & 0.51 & 0.71 $\pm$ 0.22 \\
            $\cafiftyCRNLand$      & 1.09 & 1.02 $\pm$ 0.31 & 1.24 & 1.32 $\pm$ 0.31 & 1.09 & 0.74 $\pm$ 0.31 & 1.06 & 0.76 $\pm$ 0.31 & 0.73 & 0.80 $\pm$ 0.31 \\
            $\TRPNnCRNLand$        & 2.89 & 3.04 $\pm$ 0.37 & 1.65 & 1.17 $\pm$ 0.37 & 2.14 & 1.75 $\pm$ 0.37 & 1.95 & 1.64 $\pm$ 0.37 & 2.08 & 1.75 $\pm$ 0.37 \\
            $\CVventricles$        & 0.44 & 0.64 $\pm$ 0.22 & 0.65 & 0.62 $\pm$ 0.22 & 0.55 & 0.62 $\pm$ 0.22 & 0.50 & 0.60 $\pm$ 0.22 & 0.73 & 0.84 $\pm$ 0.22 \\
            $\kFEC$                & 5.88 & 4.45 $\pm$ 1.42 & 3.07 & 3.51 $\pm$ 1.42 & 3.30 & 2.82 $\pm$ 1.42 & 5.61 & 5.95 $\pm$ 1.42 & 2.51 & 1.45 $\pm$ 1.42 \\
            $\gCaLCRN$             & 0.13 & 0.14 $\pm$ 0.05 & 0.12 & 0.11 $\pm$ 0.05 & 0.19 & 0.13 $\pm$ 0.06 & 0.14 & 0.10 $\pm$ 0.05 & 0.13 & 0.10 $\pm$ 0.05 \\
            $\btatria$             & 3.19 & 3.34 $\pm$ 0.45 & 2.40 & 2.69 $\pm$ 0.45 & 2.24 & 1.64 $\pm$ 0.45 & 2.86 & 2.28 $\pm$ 0.45 & 2.50 & 3.24 $\pm$ 0.45 \\
            $\Rsys$                & 3.28 & 3.18 $\pm$ 0.10 & 2.50 & 2.51 $\pm$ 0.10 & 3.84 & 3.76 $\pm$ 0.10 & 3.57 & 3.48 $\pm$ 0.10 & 2.26 & 2.38 $\pm$ 0.10 \\
            $\Rpulm$               & 2.63 & 2.50 $\pm$ 0.16 & 2.84 & 2.89 $\pm$ 0.16 & 3.12 & 3.07 $\pm$ 0.16 & 2.45 & 2.34 $\pm$ 0.16 & 3.40 & 3.55 $\pm$ 0.16 \\
            \bottomrule
        \end{tabular}
        \caption{$\Tall$: true value and mean plus/minus two standard deviations associated to the estimated model parameters for all the $\numTest = 5$ electromechanical simulations.}
        \label{tab:Tall}
    \end{center}
\end{sidewaystable}